\documentstyle[12pt]{amsart}

\newcommand{\ts}{\vspace{\baselineskip}\noindent{\bf Proof.}$\;\;$}

\newcommand{\cM}{{\cal M}}
\newcommand{\cO}{{\cal O}}

\newcommand{\oq}{\overline{q}}
\newcommand{\oQ}{\overline{Q}}

\newcommand{\half}{{\mbox{${1\over 2}$}}}

\newcommand{\G}{{\Gamma}}
\newcommand{\T}{{\Theta}}
\newcommand{\D}{{\Delta}}
\newcommand{\bC}{{\Bbb C}}
\newcommand{\bP}{{\Bbb P}}

\newcommand{\bT}{{\Bbb T}}
\newcommand{\bZ}{{\Bbb Z}}
\newcommand{\cI}{{\cal I}}

\newcommand{\tP}{\widetilde{\bP}}
\newcommand{\tphi}{\widetilde{\phi}}
\newcommand{\tpsi}{\widetilde{\psi}}

\newcommand{\lra}{\longrightarrow}
\newcommand{\ra}{\rightarrow}
\newcommand{\inj}{\hookrightarrow}
\newcommand{\surj}{\hspace{3pt}\to \hspace{-19pt}{\rightarrow} \:\:}

\newtheorem{proposition}{Proposition}

\newtheorem{theorem}[proposition]{Theorem}

\begin{document}
\title[Tangent space to moduli of vector bundles]{The tangent space to
the moduli space of vector bundles on a curve and the singular locus
of the theta divisor of the jacobian}
\author{B. van Geemen and E. Izadi}
\address{Dipartimento di Matematica, Universit\`a di Pavia, 
Via Ferrata 1, I-27100 Pavia, Italia}
\email{geemen@@dragon.ian.pv.cnr.it}
\address{Department of Mathematics, Boyd
Graduate Studies Research Center, University of Georgia, Athens, GA
30602-7403, USA}
\email{izadi@@math.uga.edu}

\begin{abstract} 
We complete the proof of the fact that the moduli space 
of rank two bundles with trivial determinant embeds into the linear system of divisors on 
$Pic^{g-1}C$ which are linearly equivalent to $2\Theta$.
The embedded tangent space at a semi-stable non-stable bundle $\xi\oplus\xi^{-1}$, where 
$\xi$ is a degree zero line bundle,
is shown to consist of those divisors in $|2\Theta|$
which contain $Sing(\Theta_{\xi})$ where $\T_{\xi}$ is the translate of
$\T$ by $\xi$. We also obtain geometrical results on the
structure of this tangent space.
\end{abstract}

\maketitle
\begin{center}
{\sc Introduction}
\end{center}
\vskip20pt

Let $C$ be a smooth complete irreducible algebraic curve of genus $g$,
with $g\geq 2$ over the field $\bC$ of complex numbers. For an integer
$d$, let $Pic^dC$ be the variety parametrizing line bundles of degree
$d$ on $C$. Then $Pic^{g-1}C$ has a natural theta divisor which is the
reduced divisor with underlying set
\[
\T :=\{ L \in Pic^{g-1}C:\; h^0(L) > 0\}\: .
\]
For a line bundle $L$ on $C$, let $\cM_L$ be the moduli space of
semi-stable vector bundles of rank $2$ and determinant $L$ on $C$. Put $\cO
:= \cO_C$. The variety $\cM_{\cO}$ is a normal and Cohen-Macaulay
projective variety of dimension $3g-3$ with only rational
singularities. Its singular locus is the set of (semi-stable equivalence
classes of) semi-stable non-stable bundles, i.e., bundles of the form
$\xi\oplus\xi^{
-1}$ with $\xi\in Pic^0 C$. 

The Picard group of $\cM_\cO$ is isomorphic to $\bZ$ and we let ${\cal L}$
be its ample generator. The dimensions of the vector spaces 
$H^0(\cM_\cO,{\cal L}^{\otimes k})$ are given by
the celebrated Verlinde formula.
The natural map $\cM_\cO\rightarrow \bP H^0(\cM_\cO,{\cal L})^*$
coincides with the morphism
\[
\D :\cM_{\cO}\lra | 2\T |
\]
which to a semi-stable vector bundle $E$ of rank $2$ and trivial
determinant associates the divisor $D_E$ on $Pic^{g-1}C$
defined set-theoretically as
\[
D_E :=\{ L\in Pic^{g-1} C :\; h^0(E\otimes L) > 0\}\: .
\]
We complete the proof of the following theorem:

\begin{theorem}\label{delemb}
The morphism $\D$ is an embedding for any non-hyperelliptic curve 
of genus $g\geq 4$.
\end{theorem}

Narasimhan and Ramanan showed that $\D$ is an isomorphism for
$g=2$ \cite{narasimhanramanan69} and an embedding for $g=3$ and $C$
non-hyperelliptic \cite{narasimhanramanan84}. Beauville proved
\cite{beauville88} that for $C$ hyperelliptic $\D$ induces an embedding
of the quotient of $\cM_{\cO}$ by the hyperelliptic involution of $C$
(the hyperelliptic involution induces the identity on $\cM_{\cO}$ if
$g=2$). He also proved \cite{beauville88} that if $C$ is
non-hyperelliptic, then $\D$ has generic degree $1$. Subsequently,
Laszlo proved \cite{laszlo94} that $\D$ is an embedding for $C$
general. Brivio and Verra proved \cite{brivioverra96} that $\D$ is
injective and an immersion on the smooth locus of $\cM_{\cO}$ for any
non-hyperelliptic $C$. The singular locus of $\cM_{\cO}$ can be
identified with the Kummer variety $K^0(C) := Pic^0 C/\pm 1$ via the
morphism
\[
Pic^0 C\lra\cM_{\cO },\qquad
\xi\longmapsto \xi\oplus \xi^{-1}
\]
which induces an embedding of $K^0(C)$. Using \cite{brivioverra96},
Theorem \ref{delemb} follows from Theorem \ref{emb} which we  prove in
\ref{prfgenxi} and \ref{dDinj}:

\begin{theorem}\label{emb}
Suppose $C$ is non-hyperelliptic of genus $g\geq 4$. Then, for any
$\xi\in Pic^0C$, the morphism $\D$ is an immersion at $\xi\oplus\xi^{
-1}$.
\end{theorem}

Since $\Delta$ is an embedding, we can obtain geometric information on the tangent spaces 
to $\cM_\cO$ in terms of divisors in the linear system $|2\T|$.
Let $\bT_{\xi}$ denote the embedded tangent space to $\D (\cM_{\cO })$ at
the image of $\xi\oplus\xi^{ -1}$. Then $\bT_\xi$ is a projective subspace of
$|2\Theta|$ of the same dimension as $T_{\xi\oplus\xi^{ -1}}\cM_\cO$.
Each point of $\bT_\xi$ is a divisor on $Pic^{g-1}C$. The following theorem
identifies these divisors.

\begin{theorem}\label{prel}
For 
$\xi\in Pic^0C$,
let $\T_{\xi}$ be the translate of $\T$ by $\xi$. Then
the embedded tangent space to $\Delta(\cM_\cO)$ at $\Delta(\xi\oplus\xi^{
-1})$ is:
\[
\bT_{\xi}=\{D\in |2\T|:\; Sing(\T_{\xi})\subset D\;\}.
\]
\end{theorem}

The equations for $\bT_{\xi}$ 
are linear and thus are elements of $|2\T|^*$.
The natural morphism
defined by the global sections of $2\Theta$:
\[
h:Pic^{g-1}C \longrightarrow |2\Theta|^*,\qquad L\longmapsto H_L:=
\{D\in |2\Theta|:\;L\in D\}\quad\subset |2\Theta|
\]
gives a `natural' supply of candidate equations $H_L$ for $\bT_\xi$.
Let $\langle Sing(\T_{\xi})\rangle\subset |2\Theta|^*$ be the span of
the image of $Sing(\Theta_{\xi })$ by $h$.  An equivalent form of Theorem
\ref{prel} is (cf.\ section \ref{equi}):

\begin{theorem}\label{mainthm}
For any $\xi\in Pic^0C$ the embedded tangent space $\bT_\xi$ to $\Delta(\cM_\cO)$ at $\Delta(\xi\oplus\xi^{-1})$ is the intersection of the hyperplanes $H_L$ where $L$ runs over $Sing(\T_\xi)$:
\[
\bT_{\xi} =\cap_{L\in Sing(\T_{\xi})}\, H_L\:.
\]
\end{theorem}
The difficult part of the proof is the inclusion `$\supset$' which follows
from a dimension computation in section \ref{span}.

The space $\bT_\cO$ contains much geometric information about $\cM_\cO$.
Using geometrical invariant theory
Laszlo \cite{laszlo96} proved 
\[
T_{\cO^{\oplus 2}}\cM_\cO\cong S^2H^1(\cO)\oplus \wedge^3H^1(\cO).
\]
The subspace $S^2H^1(\cO)$ of $T_{\cO^{\oplus 2}}\cM_\cO$ is the tangent space
to the Kummer variety $K^0(C)\subset\cM_\cO$. 
We define $\bT_0\;(\subset\bT)$ to be the embedded
tangent space to the image of the Kummer variety $K^0(C)$ in
$|2\Theta|$.
To investigate the $\wedge^3 H^1(\cO)$ quotient of
$T_{\cO^{\oplus 2}}\cM_\cO$
we construct a family of hyperplanes in $|2\T|$ which contain $\bT_0$
but do not contain $\bT$. These hyperplanes thus correspond to elements
in $\bP (\left(\wedge^3 H^1 (\cO )\right) )^*$.

Let $\omega$ be the canonical bundle on $C$.  There is a
morphism, similar to $\Delta$ (see \ref{candet}):
$
\Delta_\omega^*:\cM_\omega\rightarrow |2\Theta|^*
$
where $\cM_\omega$ is the moduli space of rank two bundles with canonical
determinant.
In Section $\ref{3s}$ we define a natural rational map 
$
\beta:Gr(3,H^0(\omega))\longrightarrow \cM_\omega,
$
where $Gr(3,H^0(\omega))$ is the grassmannian
of 3-planes in $H^0(\omega)$.
Composing $\beta$ with $\Delta_\omega^*$ we obtain a rational map
\[
Gr(3,H^0(\omega))\stackrel{\beta}{\longrightarrow }\cM_{\omega}
\stackrel{\Delta_{\omega}^*}{\longrightarrow} |2\Theta|^*,\qquad
W\longmapsto H_W.
\]
Note that $H^0(\omega)$ and $H^1(\cO)$ are dual vector spaces so, with
abuse of notation, $(\bT/\bT_0)^*\cong \bP\wedge^3 H^0(\omega)$. We
prove in \ref{prt2}:

\begin{theorem}\label{t2}
Let $W\in Gr(3,H^0(\omega))$ be general.
Then $\bT_0\subset H_W$. The rational restriction map:
\[
Gr(3,H^0(\omega))\longrightarrow (\bT/\bT_0)^*,
\qquad W\longmapsto \bar{H}_W:=H_W/(H_W\cap \bT_0)
\]
extends to a morphism which is the composition of the Pl\"ucker map
$W\mapsto \wedge^3W$ with a linear isomorphism $\bP (\wedge^3
H^0(\omega) )\cong
(\bT /\bT_0)^*$.
\end{theorem}

The proof of the theorem involves the construction of small
subvarieties of $\cM_\cO$ which map to points in $\bT/\bT_0$.  We
construct such subvarieties, parametrized by $p,\,q,\,r\in C$, in
section \ref{mstb} and show that we obtain a rational map $C^{(3)}\rightarrow
\bT/\bT_0=\bP \wedge^3H^1(\cO)$, which is essentially the map which
to $p+q+r$ associates the plane it spans in the canonical space. We
then study the intersection of these subvarieties with the hyperplanes
$H_W$ to prove Theorem \ref{t2} and Theorem
\ref{emb} in the case that $\xi^{\otimes 2}\cong\cO$.

Recently Pauly and Previato constructed bundles in
$\bT\cap \Delta(\cM_\cO)$ which allowed them to give another description of
the above subvarieties. They obtain a second proof of the incidence
relations in our Proposition \ref{prop411}, 
see \cite{paulypreviato}.
They also provide  generalizations of the space $\Gamma_{00}$ and
relate these to the geometry of $\cM_\cO$.

\

\noindent{\bf Acknowledgements.} We would like to thank W. M. Oxbury, C. Pauly
and E. Previato
for stimulating and helpful discussions.
The first author thanks the University of Georgia for the invitation to work
and lecture there on the results of this paper. The second author thanks the
University of Torino: it was during the second author's visit of the University
of Torino that this collaboration started.

\section{The geometry of the $|2\T|$ system.}\label{43}

Here we collect some facts on the geometry of the linear
system $|2\T|$ on $Pic^{g-1}(C)$.

\subsection{} \label{equi}
For any subset $Y$ of a
projective space $\Bbb P$, denote by $Y^{\perp }$ its polar in the
dual projective space ${\Bbb P}^*$.
Note that the linear subspace in $|2\Theta|$ defined by the 
$H_L$'s with $L\in Sing(\T_{\xi})$ is:
$$
\langle Sing(\T_{\xi})\rangle^{\perp }:=\cap_{L\in Sing(\T_{\xi})} H_L=
\{D\in |2\T|:\;Sing(\T_{\xi})\subset D\;\},
$$
this shows that Theorem \ref{mainthm} is indeed equivalent to Theorem \ref{prel}.

\subsection{Wirtinger's duality}\label{du}
The dual of the projective space $|2\Theta|$ has an intrinsic description
as a linear system on $Pic^0C$.
This is Wirtinger's duality $d_w$ (cf. \cite{mumford74} page 335):  a
canonical isomorphism
\[
d_w : |2\Theta_0|\stackrel{\cong}{\longrightarrow} |2\Theta|^*,
\]
where $\T_0$ is any symmetric theta divisor on $Pic^0 C$. It is
characterized by the property that for $\alpha\in Pic^{g-1}C$ one has
$d_w(\Theta_L+\Theta_{\omega\otimes L^{-1}})=\{D\in
|2\Theta|:\;L\in D\}$ with $\Theta_L=\{\xi\in Pic^0C:\,h^0(\xi\otimes L)>0\}$
the translate of $\Theta$ by $L$.

\subsection{Moduli of bundles with canonical determinant.}\label{candet}
Recall that $\cM_\omega$ is the moduli space of semi-stable rank 2 bundles with determinant $\omega$.
Similar to the morphism $\Delta$, we have
a morphism:
\[
\Delta_\omega:\cM_\omega\longrightarrow |2\Theta_0|,\qquad
E\longmapsto D_E :=\{L\in Pic^0C:\;h^0(E\otimes L)\neq 0\;\}.
\]
Define
\[
\D_\omega^*:\cM_\omega \lra |2\Theta|^*,\qquad
\D_\omega^*=d_w\D_\omega.
\]

Let $e_{\omega} : Pic^{g-1} C\ra\cM_{\omega}$ be the morphism which to
a line bundle of degree $g-1$ associates the vector bundle
$L\oplus(\omega\otimes L^{-1})$. Then the composition $d_w\Delta_\omega
e_\omega:Pic^{g-1}C\ra |2\T|^*$ is equal to the natural map
$h:Pic^{g-1}C\rightarrow |2\Theta|^*$.

\subsection{}
Recall that we defined $\bT_0\;(\subset\bT)$ to be the embedded
tangent space to the image of the Kummer variety $K^0(C)$ in
$|2\Theta|$. 
To recall a nice
description of the linear equations for $\bT_0$ we define $\G_{00}$ to
be the space of global sections of $\cO_{Pic^0 C} (2\T_0)$ which have
multiplicity at least $4$ at the origin:
\[
\bP\G_{ 00} :=\{ D\in |2\T_0 |: mult_\cO (D)\geq 4\}\qquad
\subset|2\Theta_0|=|2\Theta|^*.
\]
The vector space $\Gamma_{00}$ has dimension $2^g-g(g+1)/2-1$. Also define
$C-C$ to be the surface
$$
C-C:=\{\cO(p-q):\;p,\,q\in C\,\}\quad\subset Pic^0C
$$
and let $h_0:Pic^0C\rightarrow |2\Theta|=|2\Theta_0|^*$ be the natural map
associated to the linear system $|2\Theta_0|$.

\subsection{Lemma.}\label{C-Cperp}
The span of the image of $C-C\;\subset Pic^0C$ in $|2\Theta|$ is
$\bP\Gamma_{00}^{\perp}$ which is equal to $\bT_0$.

\ts Let $\delta: C\times C \lra C-C \subset Pic^0C$ be the difference
map and put $\G_0=\{s\in H^0(Pic^0C,\cO(2\T_0)):\;s(0)=0\;\}$. 
By \cite{welters86} page 18, the restriction $\delta^*
:\Gamma_0 \lra H^0(C\times C, 2\T)$ induces a surjection $\Gamma_0
\lra S^2H^0(\omega)$ which coincides with the map which to an element
of $\Gamma_0$ associates the quadratic term of its Taylor expansion at
$0$. It follows from this that the span of $C-C$ in $|2\T|$ is
contained in $\bP\Gamma_{00}^{\perp}$ and has dimension
$g(g+1)/2$. Therefore the span of $C-C$ is equal to
$\bP\Gamma_{00}^{\perp}$. 

To see that $\bP\Gamma_{00}^{\perp} = \bT_0$ or equivalently
$\bP\G_{00} = \bT_0^{\perp}$, we consider the cotangent space of the
Kummer variety at the origin. The maximal ideal of the local ring of
the Kummer variety at the origin consists of the $(-1)$-invariant
elements in the maximal ideal at the origin of $Pic^0C$.  The degree
$d$ part of the Taylor series of a regular function at $\cO$ is
canonically identified with an element of $S^dH^0(\omega)$. Therefore
the cotangent space at the origin of the Kummer variety is canonically
identified with $S^2H^0(\omega)$ and thus has dimension $g(g+1)/2$.
Since the map $r:\Gamma_0\rightarrow S^2H^0(\omega)$ is surjective,
the differential of the map $K^0(C)\rightarrow |2\Theta|$ is injective
and thus $\dim \bT_0=g(g+1)/2$. Moreover, $\bT_0$ is defined by the
kernel of $r$ which is $\bP\Gamma_{00}^{\perp}$.  \qed

\subsection{Spaces and annihilators} 
We have the following diagram of
projective subspaces and their annihilators in the dual:
\[
\begin{array}{ccccccccc}
 & &\Delta(\cO^{\oplus 2})&\subset&\bT_0&\subset&\bT&\subset&|2\Theta|\\
 &  & | & & | & & | & & \\
|2\Theta|^* & \supset & \langle \Theta \rangle & \supset & \bP
\Gamma_{00} & \supset & \langle Sing( \Theta)\rangle & &
\end{array}
\]
where $\langle\T\rangle$ and $\langle Sing(\T )\rangle$ are the spans of
the images of $\T$ and $Sing(\T)$ respectively.

\section{The general unstable bundle}\label{genunstab}

\subsection{} In this section we prove the injectivity of the differential of $\D$ at the 
non-stable bundles $\xi\oplus\xi^{-1}$ with $\xi^{\otimes 2}\not\cong\cO$. The proof is in 
Corollary \ref{prfgenxi}, the proof of 
Theorem \ref{emb} is completed in \ref{dDinj}.
We will use Bertram's results on extensions of line bundles.
These also play an important role in the case $\xi^{\otimes 2}\cong\cO$
which will be considered in section \ref{mstb}.

\subsection{Bertram's maps}\label{bertram}
For a divisor class  $D\in Pic^dC$ with $d>0$ we let
\[
C_{D}:={\rm im}\left(C\longrightarrow |\omega(2D)|^*\cong \bP H^1(\cO(-2D))
\right)
\]
be the image of $C$ under the natural map. 
Each $\epsilon\in \bP H^1(\cO(-2D))$ defines, up to isomorphism, an extension
\[
0\longrightarrow \cO(-D)\longrightarrow F_\epsilon
\longrightarrow \cO(D)\longrightarrow 0.
\]

The rational classifying map
\[
\phi_D:\bP H^1(\cO(-2D))\longrightarrow \cM_{\cO},\qquad
\epsilon\longmapsto F_\epsilon
\]
was studied in \cite{bertram92}. We denote the composition $\Delta\phi_D$ by
\[
\psi_D:
\bP H^1(\cO(-2D))\cong
|\omega(2D)|^*\stackrel{\phi_D}{\longrightarrow}
\cM_{\cO}\stackrel{\Delta}{\longrightarrow}|2\Theta| ,\qquad
\epsilon\longmapsto \Delta(F_\epsilon).
\]
 In this section we will need the case $deg(D)=1$. The case $deg(D)=2$ is studied in section \ref{bertram2}.

\subsection{The case $deg(D)=1$}\label{deg1}
In this case $\phi_D$ is an injective morphism (\cite{bertram92}, Cor.\ 4.4).
Moreover, from Theorems 1 and 2 of \cite{bertram92} one deduces that 
$\psi_D$ is given by global sections of $H^0(\bP H^1(\cO(-2D)),\cO(1))$,
hence $\psi_D$ is a linear embedding given by the complete linear system
$|\cO(1)|$. Put
\[
\bP^g_D:=\psi_D(\bP H^1(\cO(-2D)))\hookrightarrow
\Delta(\cM_{\cO})\hookrightarrow  |2\Theta|.
\]
The non-stable bundles correspond to the points of the curve
$C_D$  (\cite{bertram92}, page 451). For $p\in C_D$, we have
$\phi_D(p)=\cO(D-p)\oplus \cO(p-D)$.
This implies that $\psi_D(C_D)$ is the image of the abel-jacobi
embedded curve $C\rightarrow Pic^0C,\;p\mapsto\cO ( p-D)$ under the composition
of the maps
\[
Pic^0C\longrightarrow K^0(C)\hookrightarrow \cM_\cO
\stackrel{\Delta}{\longrightarrow} |2\Theta|.
\]
From this one concludes that the $g$-dimensional linear space $\bP^g_D$ is
the span of the image of the abel-jacobi image of $C$. See
\cite{oxburypaulypreviato} for applications of these facts to the study of
$\cM_\cO$.

\subsection{The tangent space $T_{\xi\oplus\xi^{-1}}{\cM_\cO}$}
\label{laszxi}
We consider a line bundle $\xi$ of degree zero with $\xi^{\otimes
2}\not\cong \cO$.  Let $\xi\oplus \xi^{-1}\;(\in \cM_\cO)$ be the
corresponding S-equivalence class. According to Laszlo
\cite{laszlo96} we have:
\[
T_{\xi\oplus \xi^{-1}}\cM_\cO\cong 
H^1(\cO)\oplus \left(H^1(\xi^{\otimes 2})\otimes H^1(\xi^{\otimes -2})\right),\qquad
\dim T_{\xi\oplus \xi^{-1}}\cM_\cO=g+(g-1)^2=g^2-g+1.
\]
(in fact $H^1(\cO)\cong (Ext^1(\xi,\xi)\oplus
Ext^1(\xi^{-1},\xi^{-1}))_0= (H^1(\cO)\oplus H^1(\cO))_0$ and
$H^1(\xi^{\otimes 2})=Ext^1(\xi^{-1},\xi)$). The quotient map
$Pic^0C\rightarrow K^0(C)$ induces an isomorphism from
$T_\xi Pic^0C\cong H^1(\cO)$ to the tangent space to the Kummer
variety at the image of $\xi$. The subspace $H^1(\cO)$
of $T_{\xi\oplus \xi^{-1}}\cM_\cO$ is the image of $T_\xi Pic^0C$.

\subsection{}\label{Psi}
We  define maps
$\alpha_\xi,\;\Phi_\xi$ and $\Psi_\xi$ on $C^2=C\times C$ as follows:
{\renewcommand{\arraystretch}{2.0}
\[
\begin{array}{ll}
\alpha_\xi :C\times C\lra Pic^0C ,\qquad&
(p,q)\longmapsto \xi( p-q)\; ,\\
\Phi_\xi:C\times C\stackrel{\alpha_\xi}{\lra} Pic^0C
{\lra}\cM_\cO,
\qquad&
(p,q)\longmapsto \xi( p-q)\oplus \xi^{-1}(q-p),\\
\Psi_\xi:C\times 
C\stackrel{\Phi_\xi}{\longrightarrow}\cM_\cO\stackrel{\Delta}{\longrightarrow}|2\Theta|,& 
(p,q)\longmapsto\Delta(\xi( p-q)\oplus \xi^{-1}(q-p)).
\end{array}
\]
}
Obviously, the diagonal $\D_C\subset C\times C$ is contracted to a point
by all these maps. Under $\Phi_\xi$ each  `ruling' $C\times \{q\}$
is mapped to the
curve $\phi_D(C_D)$ where $D:=\xi^{-1}(q)$. The span of
$\Psi_\xi(C\times\{q\})$ is thus the projective space
$\bP^g_D\subset \Delta(\cM_\cO)$ (see \ref{deg1}).

The following lemma relates the image of $\Psi_\xi$ to the embedded
tangent space $\bT_\xi$ to $\Delta(\cM_\cO)$ and to $Sing(\T_\xi)$.  
Next we compute the dimension of the span of $\Psi_\xi(C^2)$
by first determining the line bundle on $C^2$ which defines $\Psi_\xi$ and
its global sections (see Lemma \ref{psi*}) and then by showing that the
pull-back map is surjective on global sections
(the details of that proof are given in section \ref{2sect}).
Since, fortunately,
the dimension of this span is equal to $\dim T_{\xi\oplus\xi^{ -1}}\cM_{\cO }$,
we can deduce the injectivity of the differential in
Corollary \ref{prfgenxi}.

\subsection{Lemma.}\label{incxi}
Let $\xi\in Pic^0C$ with $\xi^{\otimes 2}\not\cong\cO$. Then we have
\[
\langle \Psi_\xi(C^2)\rangle \subset\bT_\xi\qquad{\rm and}\qquad
\langle \Psi_\xi(C^2)\rangle \subset \langle Sing(\Theta_\xi)\rangle^\perp. 
\]

\ts A point $\Psi_\xi(p,q)$ lies on the curve $\Psi_\xi(C\times\{q\})$
which lies in $\bP^g_D$ as above. 
Since $\Delta(\xi\oplus\xi^{-1})=\Psi_\xi(q,q)\in \bP^g_D$ and $\bP^g_D$ is a linear space 
in $\Delta(\cM_\cO)$ we have
$\bP^g_D\subset \bT_\xi$. Thus also $\Psi_\xi(p,q)\in\bT_\xi$.

The point $\Psi_\xi(p,q)\in|2\Theta|$ is the divisor 
$\Theta_{\xi(p-q)}+\Theta_{\xi^{-1}(q-p)}$ on $Pic^{g-1}C$.
We must show that this divisor contains $Sing(\Theta_\xi)$. For $L\in
Sing(\Theta)$, $h^0(L)\geq 2$ hence $h^0(L(q-p))\geq 1$, so
$L\in\Theta_{\cO_C(p-q)}$. Thus $Sing(\T)\subset\T_{\cO_C(p-q)}$ and so
$Sing(\Theta_\xi)\subset
\Theta_{\xi(p-q)}\subset\Theta_{\xi(p-q)}+\Theta_{\xi^{-1}(q-p)
} =\Psi_\xi(p,q)$.

This completes the proof in the case where $g\geq 5$ or $g=4$ and $C$
has two distinct pencils of degree 3. In the case where $g=4$ and $C$
has only one pencil of degree $3$, we need to be a little more careful
since $Sing\T$ is not reduced. In fact, the scheme we really need is
not $Sing\T$ but the scheme $W_{g-1}^1$ parametrizing line bundles of
degree $g-1$ with at least two sections. When $g=4$ and $C$ has only
one pencil of degree $3$, we have $Sing\T\neq W^1_{ g-1}$ whereas in
all other cases these two schemes are equal. Suppose now that $g=4$
and $C$ has only one pencil $g_3^1$ of degree $3$. The scheme $W^1_3$
has length $2$ and a one-dimensional Zariski tangent space. The
projectivization of this Zariski tangent space is the vertex of the
singular quadric containing the canonical curve. The tangent space to
$\Theta_{g^1_3(p-q)}$ at $g^1_3 $ is equal to the tangent space to
$\T$ at $\cO (g^1_3 -p + q)$. Since $C$ is non-hyperelliptic, we have
$h^0 ( g^1_3 +q ) = 2$ for every $q\in C$ and, for $p\neq q$, we have
$h^0 (g^1_3 +q -p ) = 1$. By Riemann's theorem (see for instance
\cite{ACGH} page 229), the projectivization of the tangent space to
$\T$ at $\cO (g^1_3 -p + q)$ is the span of $D_p + q$ where $D_p$ is
the unique divisor of $g^1_3$ containing $p$. This span contains the
vertex of the singular quadric containing $C$ because the quadric is
ruled by the spans of the divisors of $g^1_3$. Therefore $\T$ contains
the translate of $W_3^1$ by $\cO (p-q)$ as a scheme and
$\Theta_{\xi(p-q)}\subset\Theta_{\xi(p-q)}+\Theta_{\xi^{-1}(q-p)}$
contains the translate of $W_3^1$ by $g^1_3(p-q)$ as a scheme.  \qed

\subsection{Lemma.}\label{psi*}
Let $M:=\Psi_\xi^*\cO_{|2\Theta|}(1)$, then
\[
\dim H^0(C^2,M)=g^2-g+2.
\]

\ts The pull-back of $\cO_{|2\Theta|}(1)$ to $C^2$ under $\Psi_\xi$ 
is the pull-back of $\cO_{Pic^0 C} (2\T_0 )$ by $\alpha_{\xi }$ which is
\[
M:=\left(\pi_1^*(\omega\otimes \xi^{-2})\otimes\pi_2^*(\omega\otimes
\xi^{2})\right) (2\Delta)
\]
(where now $\Delta$ is the diagonal in $C^2$).
For the sake of completeness we sketch a proof.
For $L\in Pic^{g-1}C$ let $\T_L:=\{x\in Pic^0C:\;L\otimes x\in\T\}$
be the translate of $\T$ by $L$.
The intersection number between an abel-jacobi embedded curve and
a theta divisor is $g$, hence $\alpha_\xi(C\times\{q\})\cdot\T_L=g$.
Since $\alpha_\xi^{-1}(\T_L)=\{p\in C:\;0<h^0(L\otimes\xi(p-q))=
h^0(\omega\otimes L^{-1}\otimes\xi^{-1}(q-p))\}$,
we have 
\[
\alpha_\xi^*(\T_L)_{|C\times\{q\}}\cong \omega\otimes L^{-1}\otimes\xi^{-1}(q)
\]
for $L$ such that $h^0(\omega\otimes L^{-1}\otimes\xi^{-1}(q))=1$.
Since the pull-back is a morphism this must hold for all $L$. Thus, by
the seesaw theorem, $\alpha_\xi^*(\T_L)\cong (\pi_1^*(\omega\otimes
L^{-1}\otimes\xi^{-1})\otimes \pi_2^*N)(\Delta)$ for some line bundle
$N_\xi$ on $C$. Using the restriction to $\{p\}\times C$ one finds
similarly that $N_\xi=L\otimes\xi$. Finally, take the inverse image of
the divisor $\T_L+ \T_{\omega\otimes L^{-1}}\in|2\T|$ by $\alpha_\xi$.

In particular, we have
$M(-2\Delta)=\pi_1^*(\omega\otimes \xi^{-2})\otimes\pi_2^*(\omega\otimes
\xi^{2})$ and the K\"unneth formula gives:
\[
H^0(M(-2\Delta))\cong H^0(\omega\otimes \xi^{-2})\otimes
H^0(\omega\otimes \xi^{2})\cong \bC^{(g-1)^2},\quad
H^1(M(-2\Delta))=H^2(M(-2\Delta))=0,
\]
because $H^1(\omega\otimes \xi^{\pm 2})=0$. The exact sequence of sheaves
on $C^2$:
\[
0\lra M(-2\Delta)\lra M(-\Delta)\lra M(-\Delta)_{|\Delta}\lra 0
\]
and the isomorphisms $\Delta\cong C$, $M(-\Delta)_{|\Delta}\cong
M(-2\Delta)_{|\Delta}\otimes \cO_{C^2}(\Delta)_{|\Delta}\cong
\omega^2\otimes \omega^{-1}\cong \omega$ give:
\[
H^0(M(-\Delta))\cong H^0(\omega\otimes \xi^{-2})\otimes H^0(\omega\otimes
\xi^{2})
\oplus H^0(\omega),\qquad H^1(M(-\Delta))\cong H^1(\omega),\quad 
H^2(M(-\Delta))=0.
\]
Repeating this method we obtain:
\[
0\lra M(-\Delta)\lra M \lra M_{|\Delta}\lra 0.
\]
Since $M_{|\Delta}\cong M(-\Delta)_{|\Delta}\otimes
\cO_{C^2}(\Delta)_{|\Delta }\cong\omega\otimes\omega^{ -1}\cong\cO$, we
have $h^0(M_{|\Delta})= 1$. Since there is a section of $\cO_{|2\T|}(1)$
which is nonzero at the point $\Psi(\Delta)$, the inclusion
$H^0(M(-\Delta))\subset H^0(M)$ is strict and we conclude: $h^0(M)=1+
h^0(M(-\Delta)) = 1+ (g-1)^2 + g= g^2 - g +2$.  \qed

\subsection{Proposition.}\label{taninsingp} 
Let $\xi\in Pic^0C$ with $\xi^{\otimes 2}\not\cong\cO$
then the pull-back map
$\Psi_\xi^*:H^0(\cO_{|2\Theta|}(1))\rightarrow H^0(C^2,M)$ is surjective.
In particular
\[
\dim \langle \Psi_\xi(C^2)\rangle = g^2-g+1.
\]

\ts Restricting $\cO(1)$ to the Kummer variety and then pulling back
to $Pic^0(C)$ gives an isomorphism $H^0(\cO_{|2\Theta|}(1))\rightarrow
H^0(Pic^0C,2\T_0)$ where $\T_0$ is a symmetric theta divisor. Hence we
must show that $\alpha^*_\xi:H^0(Pic^0C,2\T_0)\rightarrow H^0(M)$ is
surjective. We will show that $\bP H^0(C^2,M)$ is spanned by the
pull-backs of the divisors $D_E$ with $E\in\cM_{\omega }$ (see
\ref{candet}).

Since the Kummer variety embeds into $|2\T|$, its embedded tangent
space has dimension $g$ at the point $p_\xi:=\Psi_\xi(\Delta )$ if
$\xi^{\otimes 2}\not\cong\cO$. The global sections of $2\Theta_0$
which are singular at $\alpha_\xi(\Delta)$ thus have codimension $g+1$
in $H^0(2\Theta_0)$ and the inverse images of these by $\alpha_{\xi }$
lie in $H^0(M(-2\Delta))$. Since $H^0(M(-2\Delta))$ also has
codimension $g+1$ in $H^0 (M)$, the map $\alpha^*_\xi$ identifies the
hyperplane sections of the embedded tangent space with $\bP
H^0(M)/H^0(M(-2\Delta))$ $\cong \bP(H^0(\omega)\oplus \bC)$
$\cong\bP^g$, hence the composition $H^0(Pic^0C,2\T_0)\rightarrow
H^0(M)\rightarrow H^0(M)/H^0(M(-2\Delta))$ is surjective.

Recall that $H^0(M(-2\Delta))\cong 
H^0(\omega\otimes \xi^2)\otimes H^0(\omega\otimes \xi^{-2})$,
and note that the projective space
$\bP (H^0(\omega\otimes \xi^2)\otimes H^0(\omega\otimes \xi^{-2}))$
is spanned by the Segre image of
$\bP (H^0(\omega\otimes \xi^2))\times\bP (H^0(\omega\otimes \xi^{-2} ))= 
|\omega\otimes \xi^2 |\times|\omega\otimes \xi^{-2} |$.
Therefore it is sufficient to prove that, for a
general element $(Z, Z')$ of
$|\omega\otimes \xi^2 |\times |\omega\otimes\xi^{-2} |$,
there is a vector bundle $E$ with $\alpha_{\xi}^{-1} (D_E ) - 2\D_C=
Z\times C + C\times Z'$. This is precisely Proposition \ref{zz'} below.
\qed

\subsection{Corollary.}\label{prfgenxi} Let $\xi\in Pic^0C$ with
$\xi^{\otimes 2}\not\cong\cO$.
Then the differential 
\[
({\rm d}\Delta)_{\xi\oplus\xi^{-1}}:T_{\xi\oplus\xi^{-1}}\cM_\cO\lra
T_{\Delta(\xi\oplus\xi^{-1})}\Delta(\cM_\cO)
\]
is an isomorphism. Thus $\dim \bT_\xi=g^2-g+1$. Moreover,
\[
\bT_\xi\subset \langle Sing(\Theta_\xi)\rangle^\perp.
\]

\ts 
The projective space $\langle \Psi_\xi(C^2)\rangle$ is spanned by the linear spaces 
$\bP^g_D$'s (with $D=\xi^{-1}(q)$ and $q\in C$). Since
$\bP^g_D=\Delta\phi_D\bP H^1(\cO(-2D))$ it lies in the image of $\Delta$,
and thus 
\[
\dim {\rm im}({\rm d}\Delta)_{\xi\oplus\xi^{-1}})\geq \dim 
\langle \Psi_\xi(C^2)\rangle=g^2-g+1=\dim T_{\xi\oplus\xi^{-1}}\cM_\cO ,
\]
the last equalities are Proposition \ref{taninsingp} and
Laszlo's result quoted in \ref{laszxi}. Hence $({\rm d}\Delta)_{\xi\oplus\xi^{-1}}$ is an 
isomorphism.

As a consequence we find the dimension of the embedded tangent space (which is equal to 
the dimension of the tangent space itself):
$
\dim \bT_\xi=g^2-g+1.
$
Since $\langle \Psi_\xi(C^2)\rangle\subset \bT_\xi$ it follows for dimension reasons that
$\langle \Psi_\xi(C^2)\rangle= \bT_\xi$.
Then  Lemma \ref{incxi} implies $\bT_\xi\subset \langle Sing(\Theta_\xi)\rangle^\perp$. 
\qed

\section{Bundles with two sections}\label{2sect}

\subsection{} In this section we finish the proof of Proposition 
\ref{taninsingp}, which follows from Proposition \ref{zz'}. We say that 
a bundle $E$ is 
generically generated by its global sections if for a general point $t$ of $C$,
there is no global section of $E$ vanishing at $t$.
Recall the map $\Psi_\xi:C^2\rightarrow |2\Theta|$, $(p,q)\mapsto
\Delta(\xi(p-q)\oplus\xi^{-1}(q-p))$. Our first result is:

\subsection{Lemma.}\label{EZZ} Suppose $E$ is a semi-stable rank 2 vector bundle with
determinant $\omega$ such that $h^0 (E\otimes\xi) = 2$ and $\alpha_\xi(C\times C)\not\subset D_E$. 

Then $E\otimes\xi$ and $E\otimes\xi^{-1}$ are
generically generated by their global sections and
\[
\alpha_{\xi}^{-1} (D_E ) - 2\D= C\times Z + Z'\times C
\]
where $Z$ and $Z'$
are the divisors on $C$ such that we have the exact sequences
\[
0\lra H^0 (E\otimes\xi )\otimes\cO\lra E\otimes\xi\lra\cO_Z\lra 0
\]
\[
0\lra H^0 (E\otimes\xi^{ -1} )\otimes\cO\lra E\otimes\xi^{ -1 }\lra\cO_{Z'
}\lra 0\; .
\]

\ts Set-theoretically $\alpha_{\xi}^{ -1} (D_E)$ is the set of pairs $(p,q)$
such that $h^0 (E\otimes\xi (p-q) ) > 0$. Therefore, if $D_E$ does not contain
the image of $\alpha_{\xi}$, then $E\otimes\xi$ is
generically generated by its global sections. Since $D_E$ is symmetric, $C- C +\xi\subset
D_E\Leftrightarrow C - C +\xi^{-1}\subset D_E$. Therefore $E\otimes\xi^{-1}$ is
also generically generated by its global sections.

It follows from the first sequence in the statement of the lemma that for
any point $t$ of $Z$ there is a section of $E\otimes\xi$ which vanishes at
$t$. In other words $h^0 ( E\otimes\xi(-t)) > 0$. Hence, for any $p\in C$,
$h^0 ( E\otimes\xi(p-t)) > 0$ and $\alpha_{\xi}^* (D_E)$ contains $C\times
Z$. Similarly, for any point $t$ of $Z'$, $h^0 ( E\otimes\xi^{ -1}(-t)) >
0$. By Riemann-Roch and because the determinant of $E$ is $\omega$, this is
equivalent to $h^0 ( E\otimes\xi( t)) > 2$. Hence $\alpha_{\xi}^* (D_E)$
contains $Z'\times C$. Since $\alpha_{\xi}^* (D_E) -2\Delta$ and $C\times Z +
Z'\times C$ are linearly equivalent, it follows that they are equal.\qed

\subsection{Lemma.}\label{nuexists}
There exists an element $(Z, Z')$ of $|\omega\otimes\xi^2
|\times |\omega\otimes\xi^{-2 }|$ and a vector bundle $E$ as in Lemma \ref{EZZ}.

\ts Since $\xi^{\otimes 2}\not\cong\cO$, there are distinct effective nonzero
divisors $D$ and $D'$ on $C$ such that $\xi^{\otimes 2}\cong\cO (D'-
D)$ and $h^0( D) = h^0 (D' ) = 1$. Since $g\geq 4$, the difference map
\[
C^{(g-2)}\times C^{(g-2)}\longrightarrow Pic^0C,\qquad
(D,D')\longmapsto\cO (D-D' )
\]
is surjective and we can find $D$ and $D'$ of degree at most $g-2$. If
$D$ and $D'$ have degree less than $g-2$, we can replace them by $D +
G$ and $D' + G$ where $G$ is a general effective divisor to obtain
deg$(D)=$ deg$(D') = g-2$, keeping the equality $h^0( D) = h^0 (D' ) =
1$. Now consider the Bertram map:
\[
\phi_{\xi^{-1}(-D)}:\bP
H^1((\omega^{-1}\otimes \xi^{\otimes
2})(2D))\longrightarrow\cM_\omega,\qquad
\epsilon\longmapsto E_\epsilon,
\]
where $E_\epsilon$ is the (S-equivalence class of) the bundle given
by the extension defined by $\epsilon$:
\[
0\longrightarrow
\xi(D)\longrightarrow E_\epsilon\longrightarrow
(\omega\otimes\xi^{-1})(-D)\longrightarrow 0.
\]
Since $deg((\omega\otimes\xi^{-1}(-D))^{\otimes 2})=2(2g-2-(g-2))=2g$,
the Bertram map is a morphism, in particular $E_\epsilon$ is indeed
semistable.

Next we tensor the defining sequence for $E_\epsilon$ by $\xi$, the
boundary map in the cohomology will be denoted by:
\[
\delta_\epsilon:H^0(\omega(-D))\longrightarrow
H^1(\xi^{\otimes 2}(D))=H^1(D')=H^0(\omega(-D'))^*.
\]
Choose $\epsilon$ such that $\delta_\epsilon$ has rank $1$. Then $h^0
(E_\epsilon\otimes\xi ) = 2 (= h^0 (E_\epsilon\otimes\xi^{ -1} ))$.

\subsubsection{Claim.}
The image of $\alpha_{\xi }$ is not contained in $D_{E_\epsilon}$
for such $\epsilon$.

\ts Set-theoretically $\alpha_{\xi}^{ -1} (D_{E_\epsilon})$ is the set
of pairs $(p,q)$ such that $h^0 (E_\epsilon\otimes\xi (p-q) ) >
0$. Let $(p,q)$ be a general element of $C\times C$. We will show that
$(p, q)$ is not contained in $\alpha_{\xi}^{ -1} (D_E)$. Since
$deg(D)=g-2$ and $h^0(D)=1$ we have $h^0(\omega(-D))=2$ and, since $q$
is general, $h^0(\omega(-D-q))=1$. Since
$deg(\omega(-D-q))=2g-2-(g-2)-1=g-1$ and $p$ is general, we then have
$h^0(\omega(-D-q+p))=1$ so that $p$ is a base point of the linear
system (of dimension zero) $|\omega(-D-q+p)|$. We have a similar
result for $D'$.

Tensoring the defining sequence of
$E_{\epsilon}$ by $\xi (p-q)$ and using $\xi^{\otimes 2}(D)\cong\cO(D')$
we obtain
\[
0\lra \cO( D' +p-q )\lra E_{\epsilon }\otimes\xi (p-q)\lra\omega (-D
+p-q )\lra 0\; .
\]
Let $\delta_{p,q} :H^0(\omega(-D+p-q))\ra H^0(\omega(-D'-p+q)^*$ be
the connecting homomorphism for the associated sequence of cohomology.
We have a diagram:
\[
\begin{array}{rcl}
H^0(\omega(-D+p-q))&\stackrel{\delta_{p,q}}{\lra}&H^0(\omega(-D'-p+q))^*\\
\cong\uparrow&&\downarrow\cong\\
H^0(\omega(-D-q))&&H^0(\omega(-D'-p))^*\\
\downarrow&&\uparrow\\
H^0(\omega(-D))&\stackrel{\delta_\epsilon}{\lra}&H^0(\omega(-D'))^*\; .
\end{array}
\]

From this one can analyze when $\delta_{p,q}$ is zero: the condition
is that $\delta_\epsilon(H^0(\omega(-D-q)))$ lies in the kernel of the map
$H^0(\omega(-D'))^*\rightarrow H^0(\omega(-D'-p)^*$. However, since
$q$ is general, the image of $H^0(\omega(-D-q))$ is a general line in
$H^0(\omega(-D))$. The image of this line by $\delta_\epsilon$ will not be
contained in the kernel of $H^0(\omega (-D') )^*\rightarrow H^0
(\omega(-D'-p)^*$ for a fixed general $p$. Therefore a general $(p,q)$
is not in $\alpha_{\xi }^{ -1} (D_{E_\epsilon} )$. \qed

\

Therefore by Lemma \ref{EZZ} we have exact sequences
\[
0\lra H^0 (E_{\epsilon}\otimes\xi )\otimes\cO\lra
E_{\epsilon}\otimes\xi\lra\cO_Z\lra 0
\]
\[
0\lra H^0 (E_{\epsilon}\otimes\xi^{-1} )\otimes\cO\lra
E_{\epsilon}\otimes\xi^{-1}\lra\cO_{ Z' }\lra 0
\]
where $(Z, Z')\in |\omega\otimes\xi^{\otimes 2}|\times
|\omega\otimes\xi^{\otimes -2}|$.
It is immediate now that for the above choice of $(Z, Z')$ and
$E = E_{\epsilon}$
the claim is true.\qed

\subsection{}\label{nuZ}
Let $Z$ be a general element of $|\omega\otimes\xi^{\otimes 2} |$. The
space of maps $\omega^{\oplus 2}\stackrel{\nu_Z }{\ra }\cO_Z$ has
dimension $2 (2g-2) = 4g-4$ and a general such map is
surjective. Consider a general such map $\nu_Z$ and let $E\otimes\xi^{
-1}$ be its kernel. Apply the functor $Hom ( .,\omega)$ to the
sequence
\[
0\lra E\otimes\xi^{ -1 }\lra\omega^{\oplus 2}\lra\cO_Z\lra 0
\]
to obtain
\[
0\lra\cO^{\oplus 2}\lra E\otimes\xi\lra\omega_{Z}\lra 0\; .
\]
As $Z$ varies in $|\omega\otimes\xi^2 |$, the spaces of surjective maps $\omega^{\oplus 
2}\stackrel{\nu_Z }{\ra }\cO_Z$ form an open subset of a vector bundle over 
$|\omega\otimes\xi^2 |$. In particular, the space of such maps is irreducible. Therefore, 
by \ref{nuexists}, for general $Z$ and a generic map $\omega^{\oplus 2}\stackrel{\nu_Z 
}{\ra }\cO_Z$, the vector bundle $E$ is semi-stable, $h^0 (E\otimes\xi ) = h^0 
(E\otimes\xi ) = 2$ and $C- C +\xi\not\subset D_E$. So, by Lemma \ref{EZZ}, we have 
$\alpha_{\xi }^{-1 }(D_E) -\Delta = C\times Z + Z'\times C$ where $Z'$ is as in Lemma 
\ref{EZZ}. Sending $\nu_Z$ to $Z'$ defines a map to $|\omega\otimes\xi^{ -2}|$ which 
identifies the space of maps $\nu_Z :\omega^{\oplus 2}\surj\cO_Z$ with the space of maps 
$\nu_{Z'} :\omega^{\oplus 2}\surj\cO_{Z'}$ since the construction is symmetric in $Z$ and 
$Z'$.

\subsection{} Let $Grass := Grass(2, H^0(\omega^{\oplus 2}))$ be the
grassmannian of two-dimensional subvector spaces of $H^0 (\omega^{\oplus
2})$. Let $U\subset Grass$ be the subset parametrizing two-dimensional
subvector spaces $V$ such that $V$ generates a rank $2$ subsheaf of
$\omega^{\oplus 2}$ and the cokernel of the evaluation map
$V\otimes\cO\lra\omega^{\oplus 2}$ is the structure sheaf of a reduced
divisor. In other words, we have the exact sequence
\[
0\lra V\otimes\cO\lra\omega^{\oplus 2}\lra\cO_D\lra 0
\]
where $D$ is a reduced divisor on $C$. It is immediate that $D\in |\omega^{\otimes 2}
|$. Sending $V$ to $D$ defines a morphism $\mu : U\ra |\omega^{\otimes 2} |$.

\subsection{Lemma.} \label{nonempty}
The image of $\mu$ is a non-empty open subset of
$|\omega^{\otimes 2} |$.

\ts Consider a divisor $D =\sum_{ i=1 }^{4g-4} p_i\in |\omega^{\otimes 2} |$
where the points $p_i$ are distinct. For any $i$, let $C_i\subset H^0
(\omega^{\oplus 2} )$ be the space of sections vanishing at
$p_i$. Requiring that the length of the cokernel of
$V\otimes\cO\ra\omega^{\oplus 2}$ at $p_i$ be at least $1$ is equivalent to
requiring $V\cap C_i\neq\{ 0\}$. Note that since the total length of the
cokernel of $V\otimes\cO\ra\omega^{\oplus 2}$ is $4g-4$, if the length at
each $p_i$ is at least $1$, then the length at $p_i$ is exactly $1$.

Let then $M_i$ be the closed subset of $Grass$ consisting of those $V$ such
that $V\cap C_i\neq\{ 0\}$. Then $D$ is in the image of $\mu$ if and only
if $(\cap_{ i=1 }^{ 4g-4 } M_i)\bigcap U\neq\emptyset$. This is clearly an
open condition on $D$ hence the image of $\mu$ is open. 
It follows from \ref{nuZ} that the image of $\mu$ is also non-empty
hence it is a dense open subset of $|\omega^{\otimes 2} |$.
\qed

\subsection{Proposition}\label{zz'}
For a general element $(Z,Z')\in |\omega\otimes \xi^2 |\times |\omega\otimes
\xi^{-2} |$, there is a semi-stable vector bundle $E$ with $\alpha_{\xi}^{-1} (D_E ) - 
2\D_C= Z\times C + C\times Z'$.

\ts
By
\ref{nuZ} the image of $\mu$ intersects the image of $| \omega\otimes\xi^2
|\times |\omega\otimes\xi^{ -2} |$. Hence by Lemma \ref{nonempty}
for a general pair $(Z,Z')$ there
is a two-dimensional vector space $V\in U$ such that we have the exact
sequence
\[
0\lra V\otimes\cO\lra\omega^{\oplus 2}\lra\cO_{Z + Z'}\lra 0\; .
\]
Define $E\otimes\xi^{ -1}$ to be the kernel of the composition
$\omega^{\oplus 2}\surj\cO_{Z + Z'}\surj\cO_Z$ so that we have the exact
sequences
\[
0\lra E\otimes\xi\lra\omega^{\oplus 2}\lra\cO_Z\lra 0
\]
\[
0\lra\cO^{\oplus 2}\lra E\otimes\xi\lra\cO_{Z'}\lra 0\; .
\]
Applying the functor $Hom ( .,\omega)$ to the first sequence, we
obtain
\[
0\lra\cO^{\oplus 2}\lra E\otimes\xi^{-1}\lra\omega_{Z}\lra 0\; .
\]
We saw that there are pairs $(Z, Z')$ for which such an $E$ exists, is
semi-stable, $h^0 (E\otimes\xi ) = 2 = h^0 (E\otimes\xi^{ -1})$ and $C-C +\xi\not\subset 
D_E$. These
are all open conditions on the space of maps $\omega^{\oplus 2}\surj\cO_Z$
(= space of maps $\omega^{\oplus 2}\surj\cO_{Z'}$) so they remain true for
a general pair $(Z, Z')$ and a general
map $\omega^{\oplus 2}\surj  \cO_{Z+Z'}$.
\qed

\section{Bundles with three sections}\label{3s}

\subsection{}\label{e3}
Bundles with 3 independent global sections are essential in
the study of $T_{\cO^2}\cM_\cO$, in particular for understanding the
$\wedge^3H^1(\cO)$ quotient of this tangent space.

Let $E$ be a vector bundle on $C$, generated by global sections,
with 
\[
{\rm rank}(E)=2,\qquad\det(E):=\wedge^2E=\omega,\qquad h^0(E)=3.
\]
For such $E$ the evaluation map
$H^0(E)\otimes_\bC\cO\rightarrow E$ is
surjective and thus we have the exact sequence:
\[
0\longrightarrow \omega^{-1}\longrightarrow H^0(E)\otimes\cO\longrightarrow
E\longrightarrow 0.
\]
Dualizing this sequence we obtain the exact sequence:
\[
0 \longrightarrow E^*\longrightarrow H^0(E)^*\otimes\cO
\stackrel{\pi}{\longrightarrow}
\omega\longrightarrow 0.
\]
The map $\pi:H^0(E)^*\otimes\cO\rightarrow \omega$ induces a 
map on global sections:
\[
\pi^0:H^0(E)^*\otimes H^0(\cO)\cong
H^0(E)^*\longrightarrow H^0(\omega);\qquad{\rm put}\quad
W_E:=Im(\pi^0).
\]

The vector space $\wedge^2H^0(E)$ is three dimensional and there is a 
natural map
\[
\alpha:\wedge^2H^0(E)\longrightarrow H^0(\wedge^2(E))=H^0(\omega).
\]
Since $H^0(E)$ is three dimensional, $H^0(E)^*\cong \wedge^2H^0(E)$.
This isomorphism can also be obtained from the Lemma below. We will
show that $\pi^0$ is injective, so $H^0(E)^*\cong_{\pi^0} W_E$, and
that $\alpha$ gives an isomorphism $\wedge^2H^0(E)\cong_{\alpha} W_E$.

\subsection{\bf Lemma.}\label{le3}
Let $E$ be a rank two bundle with $\det(E)=\omega$, $h^0(E)=3$ which is 
generated by its global sections. With the notation of \ref{e3} we have:
\begin{enumerate}
\item
The map $\pi^0$ is injective, so $\dim W_E=3$.
\item
The map $\alpha$ is injective and
\[
\alpha(\wedge^2H^0(E))=W_E.
\]
\item
For any subline bundle $L$ of $E$ we have $h^0(L)\leq 1$.
\item
The bundle $E$ is stable.
\end{enumerate}

\ts 1. The kernel of $\pi^0$ is $H^0(E^*)$. If $s\in H^0(E^*)$ is a
non-zero section, let $L$ be the subsheaf of $E^*$ generated by
$s$. Then $h^0(L)\geq 1$ and thus $deg(L)\geq 0$. The degree of
$L\otimes \omega$ is at least $2g-2$, hence $h^0(L\otimes \omega)\geq
g-1$. As $E^*\otimes \omega\cong E$, the line bundle $L\otimes
\omega$ maps into $E$ so $h^0(E)\geq g-1$. If $g>4$ this contradicts
$h^0(E)=3$. In case $g=4$ one could still have $H^0(L\otimes\omega)
=H^0(E)$, but this contradicts the fact that $E$ is globally generated.
We conclude that $H^0(E^*)=0$ and $\pi^0$ is injective.

\noindent
2.  Since $E$ is generated by global
sections, for any $p\in C$ the map on the fibers at $p$, 
$(H^0(E)\otimes \cO)_p\rightarrow E_p$, is surjective. A basis of the one-dimensional kernel
will be denoted $s_p$.

Now $s_p\in H^0(E)$ defines a two-dimensional subspace $V_p\subset 
H^0(E)^*$ which maps to a two-dimensional subspace $\pi^0(V_p)\subset 
H^0(\omega)$. It is easy to verify that 
\[
\pi^0(V_p)=\{\omega\in W_E=\pi^0H^0(E)^*:\;\omega(p)=0\,\}.
\]
Therefore the map $p\mapsto \pi^0 (V_p )$ is the composition of the
canonical map $C\ra |\omega|^*$ and the projection
$|\omega|^*\rightarrow \bP W_E^*$, with center $\bP W_E^\perp$. Since
the canonical curve spans $|\omega|^*$, the image of $C$ under this
composition spans $\bP W_E^*$.  Thus the sections $s_p$ span $H^0(E)$.

For $p,\,q\in C$ general, we have $\alpha(s_p\wedge s_q)\neq 0$ since otherwise
we would have $\alpha(s_p\wedge s_q)=0$ for all $p,\,q\in C$ which would imply that the sections $s_p$ for $p\in C$ only span a subsheaf of rank 1 of $E$, in contradiction with the fact that they span $H^0(E)$ and that $E$ is globally generated.

Let $p,\,q\in C$ be two general points. Then $V_p\cap V_q$ is one
dimensional, let $\omega_{pq}$ be a basis of $\pi^0(V_p\cap V_q)$:
\[
\langle \omega_{pq}\rangle :=\pi^0(V_p\cap V_q)=\{\omega\in 
W_E:\;\omega(p)=\omega(q)=0\,\}.
\]
Below we show that
\[
\omega_{pq}=\alpha(s_p\wedge s_q).
\]
It will follow from this that $W_E\subset \alpha(\wedge^2H^0(E))$ and,
since $\dim W_E=\dim \alpha(\wedge^2H^0(E)) = 3$ we obtain $W_E =
\alpha(\wedge^2H^0(E))$.

Let $r\in C$ be a zero of $\omega_{pq}$. Then $\omega_{pq}\in
\pi^0(V_r)$, and thus $\omega_{pq}\in V_p\cap V_q\cap V_r$. Therefore
the three sections $s_p,\,s_q$ and $s_r$ are $\bC$-linearly dependent
in $H^0(E)$. Since $s_p$ and $s_q$ are linearly independent, the
section $s_r$ is a linear combination of $s_p$ and $s_q$. So $s_p
\wedge s_q$ vanishes at $r$. In particular, the differential form
$\alpha(s_p\wedge s_q)$ is zero at $r$. We see that any zero of
$\omega_{pq}$ is also a zero of $\alpha(s_p\wedge s_q)$ and conclude
(using the general position of $p$ and $q$) that the differentials are
the same.

\noindent
3. If $L$ is a subline bundle with $h^0(L)\geq 2$, there are two 
independent sections $s,\,t\in H^0(E)$ (in the image of 
$H^0(L)\hookrightarrow H^0(E)$) whose values $s(p),\;t(p)$
are dependent in each fiber $E_p$ of $E$. Thus $s\wedge t\neq 0$ but
$\alpha(s\wedge t)=0$. This contradicts the injectivity of $\alpha$.

\noindent
4. We must show that any subline bundle $L$ of $E$ has
$deg(L)<deg(E)/2=g-1$. Suppose $L$ is a subbundle of maximal degree of
$E$ and suppose that $deg(L)\geq g-1$. The Riemann-Roch theorem
implies $h^0(L)\geq h^0(\omega\otimes L^{-1})$. The exact sequence
\[
0\longrightarrow L\longrightarrow E\longrightarrow \omega\otimes L^{-1}
\longrightarrow 0
\]
shows that $h^0(E)\leq h^0(L)+h^0(\omega\otimes L^{-1})$ thus $h^0(E)\leq 
2h^0(L)$. Since $h^0(E)=3$ we obtain $h^0(L)\geq 2$ which contradicts $3$.
\qed

\subsection{}\label{prep}
Let $W(3)\;(\subset\cM_\omega)$ be the locus of stable bundles $E$ of rank two
on $C$ with $det(E)=\omega$ and $h^0(E)=3$. We denote by
$W(3)_+\;(\subset W(3))$ the locus of bundles in $W(3)$ which are generated
by global sections.

To determine $W(3)_+$ we use the injectivity of
$\pi^0:H^0(E)^*\hookrightarrow H^0(\omega)$. For any three-dimensional
subspace $W\subset H^0(\omega)$ we can define a rank two bundle
$E_W^*$ as
\[
E_W^*:=\ker(\pi:W\otimes \cO\longrightarrow \omega)
\]
and let $E_W:=(E_W^*)^*$ be the dual bundle (here $\pi$ is the
evaluation map). The determinant of $E_W$ is $\omega$ if and only if
$\pi$ is surjective. If $\pi$ is surjective, we have an exact
sequence: $$ 0\longrightarrow \omega^{-1}\longrightarrow
W^*\otimes\cO\longrightarrow E_W\longrightarrow 0\qquad\qquad(*) $$ so
$E_W$ is generated by global sections.

\subsection{}
Denote the grassmannian of three-dimensional subspaces of 
$H^0(\omega)$ by $Gr(3,H^0(\omega))$. It has dimension $3(g-3)$ and its
Picard group is generated by an ample line bundle which we denote as $\cO(1)$. 
The Pl\"ucker embedding is the natural map 
\[
Gr(3,H^0(\omega))\longrightarrow \bP H^0(Gr(3,H^0(\omega)),\cO(1))^*
\cong \bP \wedge^3H^0(\omega).
\]
Let
\[
B_C\;\subset Gr(3,H^0(\omega))
\]
be the locus of $W\in Gr(3,H^0(\omega))$ such that the two-dimensional
linear system $|W|\;(\subset |\omega|)$ has a base point. We will
denote by
\[
D_C\;\subset Gr(3,H^0(\omega))
\]
the locus of $W\in Gr(3,H^0(\omega))$ such that the multiplication map
\[
m_W:W\otimes H^0(\omega)\longrightarrow H^0(\omega^{\otimes 2})
\]
is not surjective.
Note that we have $B_C\subset D_C$ since $W\in B_C$ implies that the quadratic
differentials in the image of $m_W$ all have a zero at a base point of
$|W|$.

\subsection{\bf Proposition.}\label{win} 
\begin{enumerate}
\item
The map
\[
S:W(3)_+\longrightarrow Gr(3,H^0(\omega))-D_C,\qquad
E\longmapsto W_E
\]
is an isomorphism, its inverse is the map
\[
\beta:Gr(3,H^0(\omega))-D_C\longrightarrow \cM_\omega,\qquad W\longmapsto E_W.
\]
\item
The locus $D_C$ has codimension one in $Gr(3,H^0(\omega))$ and $D_C$ is the 
divisor of a global section of $\cO(g-2)$.
\end{enumerate}

\ts
1. First we show that $Im(S)\subset Gr(3,H^0(\omega))-D_C$. 
From the cohomology of the  exact sequence $(*)$ in \ref{prep}
we see that $h^0(E)=3$ implies that the map
$$
H^1(\omega^{-1})\longrightarrow W_E^*\otimes H^1(\cO)
$$
is injective. Thus the dual of this map, which is the multiplication map
$m_{W_E}$,
is surjective.

From the constructions it is clear that, if $\beta$ is well-defined,
the morphisms $S$ and $\beta$ are each other's inverses. To prove that
$\beta$ is well-defined, consider an element $W$ of
$Gr(3,H^0(\omega))-D_C$. Since the multiplication map $m_W$ is
surjective, for each $p\in C$ there is an $\omega\in W$ with
$\omega(p)\neq 0$. Then the evaluation map $e_W : W \otimes \cO
\longrightarrow \omega$ is surjective so $E_W$ has determinant
$\omega$ and is generated by global sections. The cohomology of $(*)$
and the fact that $m_W$ is surjective show that $h^0(E_W)=3$. From
Lemma \ref{le3} we see that $E_W$ is stable.

2. Lazarsfeld (Theorem 1.1 of \cite{gieseker92}) and Beauville (private
communication) proved that for any non-hyperelliptic
curve there is a three-dimensional subspace $W$ of $H^0(\omega)$ such
that the multiplication map $m_W$ is surjective.
Thus $D_C$ is a Zariski closed subset of codimension $\geq 1$.

For non-hyperelliptic curves the multiplication map: $
m:S^2H^0(\omega)\longrightarrow H^0(\omega^{\otimes 2}) $ is
surjective, so its kernel, which we denote by $I_2$, has dimension
$\half(g-2)(g-3)$. Let $S_W:=Im(W\otimes H^0(\omega)\rightarrow
S^2H^0(\omega))$. Then $D_C$ is the locus of $W \in Gr(3,
H^0(\omega))$ such that $I_2 \cap S_W \neq \emptyset$. Since the
kernel of $W\otimes H^0(\omega)\rightarrow S^2H^0(\omega)$ is $(W
\otimes H^0(\omega)) \cap \wedge^2H^0(\omega) = \wedge^2 W$, we have
$\dim S_W = 3g-3$. Therefore $I_2$ and $S_W$ have complementary
dimensions and $D_C$ is the divisor of zeros of the pull-back of a
section of the Pl\"ucker bundle on $Gr(3g-3,S^2H^0(\omega))$ under the
map:
\[
Gr(3,H^0(\omega))\longrightarrow Gr(3g-3,S^2H^0(\omega)),\qquad
W\longmapsto S_W=(W\otimes H^0(\omega))/\wedge^2W.
\]
Recall that the Pl\"ucker bundle is the determinant of the universal
quotient bundle over the Grassmanian. It is also the dual of the
determinant of the universal subbundle. Let ${\bf W}$ be the universal
subbundle on $Gr(3,H^0(\omega))$. The vector bundle ${\bf
W}\otimes H^0(\omega)$ (where we view $H^0(\omega)$ as a trivial
bundle of rank $g$ over the Grassmanian) has determinant $det({\bf
W})^g=\cO(-g)$. Define a bundle ${\bf S}$ on $Gr(3,H^0(\omega))$ by the
exact sequence:
\[
0\longrightarrow \wedge^2{\bf W}\longrightarrow {\bf W}\otimes H^0(\omega)
\longrightarrow {\bf S}\longrightarrow 0.
\]
Then the fiber of ${\bf S}$ over $W\in Gr(3,H^0(\omega))$ is the
vector space $S_W$ and ${\bf S}$ is the pull-back of the universal
subbundle on $Gr(3g-3,S^2H^0(\omega))$. The dual of the determinant of
${\bf S}$ is $det({\bf W}\otimes H^0(\omega))^{-1}\otimes
det(\wedge^2{\bf W})$. For a bundle $F$ of rank $k$, the determinant
of $\wedge^2 F$ is $det(F)^{\otimes (k-1)}$ (use the splitting
principle for example). Thus $det(\wedge^2{\bf W})=det({\bf
W})^{\otimes 2}=\cO(-2)$ and $det({\bf S})^{-1}=\cO(g-2)$.
\qed

\subsection{\bf Remark.} Let $W_3\subset \cM_\omega$ be locus of all
stable bundles $E$ with $h^0(E)\geq 3$. From the results of
\cite{collinopirola95} one has $\dim W_3=3(g-3)\;(=\dim
Gr(3,H^0(\omega)))$. In fact their results imply that the
closure of $W(3)_+$ is the
only component of $W_3$ of maximal dimension.

\subsection{\bf Examples.}\label{exa}
\begin{enumerate} \item In case $g=3$, $Gr(3,H^0(\omega))$ is a point and
$S^2H^0(\omega)\longrightarrow H^0(\omega^{\otimes 2})$ is an
isomorphism. The corresponding bundle $E$ was discovered by Laszlo
\cite{laszlo91} who also showed that $W(3)=\{E\}$, thus for a
(non-hyperelliptic) genus three curve $W(3)=W(3)_+$.

\item In case $g=4$, $Gr(3,H^0(\omega))\cong |\omega|^*,\;W\mapsto
W^\perp$. The canonical curve $C_{can}$ lies on a unique quadric $Q$
which is thus $D_C$ and $B_C=C_{can}$. From \cite{oxburypaulypreviato}
we know that $\beta$ is the rational map
\[
\beta:|\omega|^*\longrightarrow |2\Theta|^*,\qquad
x=(X_0:\ldots:X_4)\longmapsto
(X_0Q(x):\ldots:X_3Q(x):R(x))\quad(\subset\bP^4\subset |2\Theta|^*)
\]
where $R$ is any cubic such that $C_{can}$ is defined by $Q$ and $R$.
The image of $\beta$ is a cubic threefold $X$ with one node in $\bP^4$.
The inverse of $\beta$ is projection from the node $X\rightarrow 
\bP^3\cong |\omega|^*$.

In case $W^\perp\in Q \setminus C_{can}$, the bundle $E_W$ is in the
S-equivalence class of the semi-stable non-stable bundle $g^1_3\oplus
h^1_3$ where $\{g^1_3,\;h^1_3\}$ is the set of line bundles $L$ of
degree $3$ such that $h^0(L)=2$. These bundles $E_W$ map to the node of $X$.

In case $W^\perp=p\in C_{can}$, the bundle $E_W$ has determinant $\omega(-p)$.
The Hecke transforms of $E_W$ at $p$, tensored by $\cO(p)$, have 
determinant $\omega$. In this way each point of $C$ determines a line in 
$X\subset \cM_\omega$. These lines all pass through the node of $X$
and under the projection 
from the node they are contracted to the corresponding points on $C$.

\item In case $g=6$ and $C$ is general, the locus $W_3$ has 6
irreducible components, one is the closure of the image of
$\beta$. The other 5 are threefolds (cones over $\bP^2$), one for
each of the 5 $g^1_4$'s on the curve, any bundle in such a component
has a subbundle (the $g^1_4$) with 2 sections (\cite{oxburypaulypreviato}).
\end{enumerate}

\subsection{} The following lemma will be used in the proof of Theorem
\ref{t2}.

\subsection{Lemma.}\label{secEt}
Let $E\in W(3)_+$ and let $p,\;q\in C$. Then:
\[
h^0(E(p))=3\;\Longleftrightarrow \; \dim(W_E\cap H^0(\omega(-p)))=2.
\]
Moreover, $h^0(E(p+q))\geq 4$ and
\[
h^0(E(p+q))=4\;\Longleftrightarrow \; \dim(W_E\cap H^0(\omega(-p-q)))=1.
\]
The intersections take place in $H^0(\omega)$. 

Furthermore, if $h^0(E(p+q))=4$, then $h^0(E(p))=h^0(E(q))=3$ and the bundle
$E(p+q)$ has a global section which is non-zero in $p$ and in $q$.

\ts
Let $D=p$. Since $\chi(E(D))=deg(E(D))+2(1-g)=2g+2(1-g)=2$, we have
$h^0(E(D))=2+h^1(E(D))$. So $h^0(E(D)) = 3 \Leftrightarrow
h^1(E(D))=1$. We have the exact sequence
\[
0\longrightarrow \omega^{-1}(D)\longrightarrow
H^0(E)\otimes\cO(D)\longrightarrow E(D)\longrightarrow 0
\]
with cohomology sequence
\[
0\longrightarrow H^0(E)\longrightarrow H^0(E(D))\longrightarrow
H^1(\omega^{-1}(D))\longrightarrow H^0(E)\otimes H^1(D)\longrightarrow
H^1(E(D)) \longrightarrow 0.
\]
Dualizing part of the cohomology sequence
and using $H^0(E)^*\cong W_E\;(\subset H^0(\omega))$ we see:
\[
H^1(E(D))^*=\ker\left(m_D:W_E\otimes H^0(\omega(-D))\longrightarrow 
H^0(\omega^{\otimes 2}(-D))\right),
\]
where $m_D$ is the multiplication map. Using a nonzero global section
of $\cO(D)$ we obtain inclusions $H^0(\omega(-D))\subset H^0(\omega)$
and $H^0(\omega^{\otimes 2}(-D))\subset H^0(\omega^{\otimes 2})$. The
map $W_E\otimes H^0(\omega)\lra H^0(\omega^{\otimes 2})$ is surjective
with kernel $\wedge^2 W_E$. Since $\dim W_E=3$ and
$h^0(\omega(-D))=g-1$ we obtain $\dim W_E\cap H^0(\omega(-D))\geq 2$.  If
$\dim (W_E\cap H^0(\omega(-D)))>2$, we have $W_E\subset
H^0(\omega(-D))$ and thus $\wedge^2W_E\subset \ker(m_D)$. So
$h^1(E(D)) \geq 3$. If $\dim(W_E\cap H^0(\omega(-D)))=2$, and
$s,\,t\in W_E\cap H^0((\omega(-D))$ are independent, we have $\langle
s\wedge t\rangle=(\ker(m_D))=1$ and thus $h^1(E(D))=1$.

Now let $D=p+q$, since $C$ is not hyperelliptic, $h^0(D)=1$ and the
proof is similar: Since $\chi(E(D))=2g+2+2(1-g)=4$, we have
$h^0(E(D))=4+h^1(E(D))$. Now $\dim W_E=3,\;\dim H^0(\omega(-D))=g-2$
so $\dim W_E\cap H^0((\omega(-D))\geq 1$. If $\dim W_E\cap
H^0(\omega(-D))>1$, and $s,\,t\in W_E\cap H^0((\omega(-D))$ are
independent, we have $s\wedge t\in \ker(m_D)$ and thus
$h^0(E(D))>4$. If $\dim W_E\cap H^0(\omega(-D))=1$ the map $m_D$ is
injective and thus $h^0(E(D))=4$.

Finally we observe that if $h^0(E(p+q))=4$ we have $\dim W_E\cap
H^0(\omega(-p-q))=1$ and thus $\dim W_E\cap H^0(\omega(-p))=2$, so
$h^0(E(p))=3$. Similarly, $h^0(E(q))=3$. Thus there is a $t\in
H^0(E(p+q))$ with $t(p)\neq 0$ and $t(q)\neq 0$.
\qed

\section{Bertram's map in degree two}\label{bertram2}

\subsection{} In this section we study Bertram's extension map
in the case $deg(D)=2$ and $D$ effective. The main result is Proposition
\ref{propbertram2} which we will use in the next section.
Write $D=p+q$ with $p,\;q\in C$. 
We recall some results of \cite{bertram92}, using the notation of
section \ref{genunstab}. 

Let $\epsilon\in H^1(\cO(-2D))$. An extension $F_\epsilon$ is stable
if $\epsilon$ does not lie on any secant line of $C_D$.  The extension
is semi-stable, but not stable, if $\epsilon$ lies on a secant line
$\langle r,s\rangle$ of $C_D$ but $\epsilon\not \in C_D$.  In this
case $F_\epsilon$ is Sheshadri equivalent to
$\cO(p+q-r-s)\oplus\cO(r+s-p-q)$. In case $\epsilon\in C_D$,
$F_\epsilon$ is unstable and $\phi_D$ is not defined on the curve
$C_D$.

From \cite{bertram92}, Theorem 1 (and p.\ 451-452), 
we know that $\phi_D$ lifts to a morphism
\[
\tphi_D :\tP^{g+2 }\lra\cM_{\cO }
\]
where $\tP^{g+2}$ is the blow-up of $\bP H^1(\cO(-2D))$ along $C_D$.
Composing the map $\tphi_D$ with $\Delta$ gives a morphism:
\[
\tpsi_D :=\D\tphi_D : \tP^{g+2}
\longrightarrow |2\Theta|.
\]
Let
\[
\bP_D:=\langle \tpsi_D(\tP^{g+2 })\rangle
\]
be the span of the image of this morphism. 
The fiber over $r\in C$ of the blow-up morphism $\tP^{g+2}
\rightarrow \bP^{g+2}$ is $|\omega(2r)|^*$ and the image in $|2\T|$ of 
$|\omega(2r)|^*$ 
is the linear space $\bP^g_r$ (as in section \ref{deg1}).

\subsection{} \label{psisym}
To determine the coordinate functions of $\tpsi_D$ we use
the following observation.
As we saw above, given two points $r,\,s\in C$, the secant line
$<r,s>\in\bP H^1(\cO(-2D))$ is mapped to the bundle
$\cO(r+s-D)\oplus\cO(D-r-s)$. Thus the image of
$\tpsi_D$ contains the image of the surface $C^{(2)}$
(symmetric product) under the composition of the abel-jacobi map
\[
\alpha_D:C^{(2)}\lra Pic^0(C),\qquad r+s\longmapsto r+s-D
\]
with the map $Pic^0C\rightarrow K^0(C)\subset \cM_\cO$ and finally
$\Delta:\cM_\cO\rightarrow |2\T|$. We simply write 
$$
\Delta\alpha_D:C^{(2)}\longrightarrow |2\Theta|,\qquad
r+s\longmapsto \T_{D-(r+s)}+\T_{(r+s)-D}
$$
for this composition.  Let $M_D$ be the pull-back of
$\cO_{|2\T|}(1)$ to $C^{(2)}$, then $M_D$ is also
$\alpha_D^*(\cO(2\Theta_0))$ where $\Theta_0$ is a symmetric theta
divisor on $Pic^0C$.  We are going to show that
$\Delta\alpha_D(C^{(2)})$ spans $\bP_D$ and this, combined with
results of Bertram, allows us to determine $\tpsi_D$ in Proposition
\ref{propbertram2}. Lemma \ref{two} and
Proposition \ref{pulsurj} below are due to C. Pauly.

\subsection{Lemma.}\label{dimmd} For any $D\in C^{(2)}$ we have:
$$
\dim H^0(C^{(2)},M_D)=1+g(g+1)/2.
$$

\ts
See \cite{brivioverra96}, Prop.\ 4.9 (and also \cite{oxburypauly}, Prop. 10.1).
\qed

\subsection{Lemma. (C. Pauly)}\label{two}
Let $p_0,\ldots,p_g$ be $g+1$ general points in $C$. Then
the images of the $p_i+p_j$, $0\leq i<j\leq g$ span a hyperplane in the space
$|M_D|^*$.

\ts Let $H_{ij}\;(\subset |M_D|)$ be the hyperplane defined by the
point $p_i+p_j$. Proving the lemma is equivalent to showing that $
\cap_{i<j} H_{ij}$ is one point. For dimension reasons, it is at least
one point. So we need to show it is at most one point.

We use the canonical isomorphism from \cite{oxburypauly}, prop.\ 10.2:
$$
\begin{array}{ccc}
H^0(C^{ (2)},M_D)&\stackrel{\cong}{\lra}&I_C(2):=H^0(|\omega(2D)|^*, {\cal
I}_C(2))\cr
\end{array}
$$
where ${\cal I}_C$ is the ideal sheaf of $C_D$ in in $|\omega(2D)|^*$.
The hyperplane $H_{ij}\;(\subset |M_D|)$ corresponds to the hyperplane
of quadrics in $|\omega(2D)|^*$ containing the curve and the secant
line $<p_i,p_j>$. Therefore $ \cap_{i<j} H_{ij} $ consists of the
quadrics in $|\omega(2D)|^*$ which contain $C_D$ and the span of
$p_0,\ldots,p_g$. We claim there is at most one such quadric:

Consider a general hyperplane $H\subset |\omega(2D
)|^*\;(\cong\bP^{g+2})$ containing the span of $p_0,\ldots,p_g$ and
let $q_1,\ldots,q_{g+1}$ be the residual points of intersection of $H$
with $C$. Since the curve spans $\bP^{g+2}$, any quadric $Q\in I_C(2)$ is
irreducible and so $H$ cannot be a component of $Q$. Thus the
restriction map $I_C(2)\rightarrow H^0(H,\cO_H(2))$ is injective. Any
quadric containing the curve and the span of the $p_i$'s intersects
$H$ in a quadric which contains the $q_i$'s and the span of the
$p_i$'s. Hence it is enough to show that there is only one such
quadric $Q$ in $H$. Since each secant line $<q_i,q_j>$ is contained in
$H$, it intersects the span of the $p_i$. Hence $<q_i,q_j>$ meets $Q$
in at least 3 points and is thus contained in $Q$. Since $Q$ contains
all the $<q_i,q_j>$ it also contains the span of the $q_i$ which is a
$\bP^g$. Thus $Q\cap H\supset \langle p_0,\ldots,p_g\rangle\cup\langle
q_1,\ldots,q_{g+1}\rangle$ and we have equality since the two sides
have the same degree and dimension and the right-hand side is reduced.
\qed

\subsection{Proposition. (C. Pauly)}\label{pulsurj}
Let $D\in Pic^2(C)$ and let
$$
\alpha_D:C^{ (2)}\lra Pic^0C,\qquad r+s\longmapsto r+s-D.
$$
Then the pull-back map
$$
\alpha_D^*:H^0(Pic^0C,\cO(2\Theta_0))\lra H^0(C^{ (2)},M_D)
$$
is surjective. Hence $\dim\langle\Delta\alpha_D (C^{(2)} )\rangle = g(g+1)/2$.

\ts By Lemma \ref{two} and with the notation there, the intersection
of the $g(g+1)/2$ hyperplanes $H_{ij}$ in $|M_D|^*$ is one point. This
intersection is
$\alpha_D^*(\Theta_\xi+\Theta_{\omega\otimes\xi^{-1}})$, with
$\xi:=\cO(p_0+\ldots+p_g-D)\in Pic^{g-1} C$ because:
\begin{enumerate}
\item $\alpha_D ( C^{ (2) })\not\subset\T_{\xi }
+\T_{\omega\otimes\xi^{ -1}}$ since for general $r, s\in C$: $h^0 (\xi
(r+s - D)) = h^0 (p_0 +\ldots + p_g +r+s - 2D) = 0$ (the $p_i$'s are
also general) and $h^0 (\omega\otimes\xi^{ -1 } (r+s - D)) = h^0 (\xi
(D -r-s )) = h^0 ( p_0 +\ldots +p_g - r- s) = 0$.
\item for all $i<j$, we have $h^0(\xi(D-p_i-p_j))>0$, hence $h^0
(\omega\otimes\xi^{ -1 } (p_i +p_j - D)) > 0$, hence $\cO (p_i +p_j -D
)\in\Theta_{\omega\otimes\xi^{ -1}}$
which implies $\alpha_D^*(\Theta_\xi+\Theta_{\omega\otimes\xi^{-1}})\in
H_{ij}$.
\end{enumerate}

Next we construct bundles $E_{ij}\in \cM_\omega$ with $0\leq i<j\leq
g$ whose divisors $D_{ij}:=\Delta(E_{ij}) = D_{ E_{ ij}}\;(\subset
Pic^0C)$ (see \ref{candet}) have the property:
$$
\cO (p_i+p_j-D )\not\in D_{ij},\qquad \cO (p_k+p_l-D )\in D_{ij}\quad{\rm
if}\;\{k,l\}\neq\{i,j\}.
$$
These conditions are equivalent to:
$$
h^0(E_{ij}(p_i+p_j - D))=0,\qquad h^0(E_{ij}(p_k +p_l - D))>0
\quad{\rm if}\;\{i,j\}\neq\{k,l\},
$$
or, equivalently, using Serre Duality, Riemann-Roch and $det(E_{
ij})\cong\omega$,
$$
h^0(E_{ij}(D-p_i-p_j))=0,\qquad h^0(E_{ij}(D-p_k-p_l))>0
\quad{\rm if}\;\{i,j\}\neq\{k,l\}.
$$
For simplicity we take the indices to be $i=0,\,j=1$. 
We consider extensions:
\[
0\lra y\lra E\lra \omega\otimes y^{-1}\lra 0,\qquad
y=\cO(p_1+\ldots +p_g-D)\;(\in Pic^{g-2}C).
\]
Then $E$ is stable since $deg(\omega\otimes y^{-2})=2$
(cf. \ref{deg1}).
From the exact sequence it is clear that
\[
h^0(E(D-p_i-p_j))\geq
h^0(p_1+\ldots+\hat{p}_i+\ldots+\hat{p}_j+\ldots+p_g)>0 \qquad{\rm
  for}\;1\leq i<j\leq g,
\]
hence $\cO (p_i+p_j-D)\in \D(E)$ if $1\leq i<j\leq g$.

It remains to consider the points $p_0+p_i$. Put
$\lambda_i=\cO(D-p_0-p_i)$ for $1\leq i\leq g$. Then, for general
$p_i$'s, $h^0(y\otimes\lambda_i)=0$ and $h^0(\omega\otimes
y^{-1}\otimes\lambda_i)=1$. The long exact sequence
\[
0\lra H^0(y\otimes\lambda_i)\lra H^0(E\otimes \lambda_i)
\lra H^0(\omega\otimes y^{-1}\otimes\lambda_i)
\stackrel{\delta(\epsilon)}{\lra}H^1(y\otimes\lambda_i)\ldots
\]
shows that $h^0(E\otimes\lambda_i)>0$ iff $\delta(\epsilon)=0$
where $\epsilon \in H^0(\omega^2\otimes y^{ -2})^*$ is the extension
class defining $E$. This coboundary is zero exactly when the image of
the multiplication map
$$
m_i:H^0(\omega\otimes y^{-1}\otimes\lambda_i)\otimes H^0(\omega\otimes 
y^{-1}\otimes\lambda_i^{-1})\lra H^0(\omega^2\otimes y^{-2})
$$
is contained in the hyperplane $H_\epsilon\subset H^0 (\omega^2\otimes
y^{ -2})$ defined by $\epsilon$.

\subsubsection{Claim.} For general $p_i$'s 
the images of the $m_i$'s are independent.

\ts We have $\omega\otimes y^{-1}\otimes\lambda_i =\omega ( 2 D - p_0
-\ldots - 2 p_i -\ldots -p_g )$ and $\omega\otimes
y^{-1}\otimes\lambda_i^{ -1} =\omega (p_0 - p_1\ldots - \hat{p_i}
-\ldots -p_g )$. Since the $p_i's$ are general, all these sheaves have
exactly one non-zero global section. For $1\leq i\leq g$, let $s_i$ be
a section of $\omega$ vanishing at $p_j$ for $1\leq j\leq g, j\neq i$
and let $t_i$ be a section of $\omega (2D -p_0\ldots - p_g )$
vanishing at $p_i$. Then, since the $p_i$'s are general, the sections
$s_i$ generate $H^0 ( C,\omega )$ and the sections $t_i$ generate $H^0
(C,\omega (2D -p_0\ldots - p_g ))$. By the
base-point-free-pencil-trick, the multiplication map $H^0 ( C,\omega
)\otimes H^0 (C,\omega (2D -p_0\ldots - p_g ) )\ra H^0 ( C,\omega^2
(2D -p_0\ldots - p_g ))$ is injective. Since the $s_i\otimes t_i$ are
independent, so are their images $s_i t_i\in H^0 ( C,\omega^2 (2D
-p_0\ldots - p_g ))$. Now, since $p_0$ is the base point of $|\omega
(p_0 )|$, the image of $H^0(\omega\otimes
y^{-1}\otimes\lambda_i)\otimes H^0(\omega\otimes
y^{-1}\otimes\lambda_i^{-1})$ in $H^0(\omega^2\otimes y^{-2})$ is
generated by $s_i t_i s_0$ where $s_0$ is a non-zero section of $\cO
(p_0 )$. Since the $s_it_i$ are independent, so are the $s_it_is_0$.
\qed
\vskip10pt

Thus we can find an $\epsilon$ with
\[
{\rm im}(m_1)\not\in H_\epsilon,\qquad {\rm im}(m_i)\in H_\epsilon
\quad {\rm for}\quad i\geq 2.
\]
The bundle $E$ defined by such an $\epsilon$ has the properties
required of $E_{01}$.

We conclude that the family 
$$
\{\alpha^*_D(\Theta_\xi+\Theta_{\omega\otimes\xi^{-1}})\}\cup
\{\alpha^*_D(D_{ij})\}_{0\leq i<j\leq g}
$$
generates $H^0(C^{ (2)},M_D)$, hence $\alpha^*_D$ is surjective.
\qed

\subsection{Corollary.}\label{propbertram2}
The composition (a rational map)
\[
\psi_D :=\D\phi_D :\bP H^1(\cO(-2D))
\longrightarrow |2\Theta|
\]
is given by the linear system of all quadrics containing
$C_D\subset |\omega(2D)|^*$.
The dimension of $\bP_D$ is $\half g(g+1)$.

\ts
From Theorems 1 and 2 of \cite{bertram92} one deduces that 
the rational map $\psi_D$ is given by a subspace of 
$I_C(2)=H^0(\bP H^1(\cO(-2D)),\cI_C(2))$ (notation as in the proof of
Lemma \ref{two}). 
In particular, $\dim
\bP_D\leq \dim |I_C(2)|$ and we have equality only if $\psi_D$ is given
by all quadrics in $I_C(2)$.
Since $I_C(2)$ is the kernel of the 
(surjective) multiplication map
$$
S^2H^0(\omega(2D) )\longrightarrow H^0(\omega^{\otimes 2}(4D)),
$$
we get from Riemann-Roch that $\dim |I_C(2)|=g(g+1)/2$.
Since $\D\alpha_D (C^{(2)} )$ lies in $\bP_D$,
we get $\dim \bP_D\geq g(g+1)/2$ from \ref{pulsurj}.
The Proposition follows.
\qed

\section{The moduli space near the trivial bundle}\label{mstb}

\subsection{}
We use Bertram's extension maps (see \ref{bertram}), restricted to
certain three-dimensional projective spaces, to study $\cM_\cO$ near
the (very) singular point $\cO^2$. The image of such a $\bP^3$ in
$|2\T|$ spans a $\bP^4$ which contains $\Delta(\cO^2)$ and which lies
in $\bT$. This $\bP^4$ intersects $\bT_0$ in a three-dimensional
space and thus defines a point in $\bT/\bT_0$. An investigation of the
intersection of such a $\bP^4$ with the hyperplanes $H_W$, with $W\in
Gr(3,H^0(\omega))$ leads to the proofs of Theorems
\ref{emb} and \ref{t2}.

\subsection{}\label{notpqr}
Choose a divisor $D=p+q\;(\;\in C^{(2)})$. Since $C$ is not
hyperelliptic, we have $h^0(p+q)=1$. Using the square of a non-zero
section $s\in H^0(\cO(D))$ we obtain a natural inclusion
\[
H^0(\omega)\hookrightarrow H^0(\omega(2D))\; .
\]
Let
\[
\pi_D:H^1(\cO(-2D))\surj H^1(\cO),
\]
be the dual map. Then the kernel of $\pi_D$ has dimension 3. Thus $\bP
\ker(\pi_D)\cong \bP^2\hookrightarrow |\omega(2D)|^*$. By Proposition
\ref{theta} below, $\bP \ker(\pi_D)=\langle 2p+2q\rangle$. Let $U'$ be
the open subset of $C^{ (2)}\times C$ parametrizing points $(p+q, r)$
such that $p,q,r$ are distinct and
$r\not\in\langle 2p+2q\rangle$. Then, for $(p+q, r)\in U'$,
\[
\langle 2p+2q+r\rangle=\langle r,\;\bP \ker(\pi_D)\rangle
\qquad\subset |\omega(2D)|^*.
\]
Proposition \ref{theta} determines
the restriction of $\psi_D$ to the three-dimensional projective space
$\langle 2p+2q+r\rangle$.

\subsection{Proposition.} \label{theta} Suppose $(p+q, r)\in U'$. We have:
\begin{enumerate}
\item The projective plane $\bP \ker(\pi_D)$ is spanned by the secant
$l:=\langle p,\;q\rangle$ and the tangent lines $l_p,\;l_q$ at $p$ and $q$
to $C_D\;(\subset |\omega(2D)|^*)$.

\item
The map $\psi_D$ restricted to $\bP \ker(\pi_D)$ is given by the pencil 
of conics $l^2$ and $l_pl_q$. 

\item
The image under $\tpsi_D$ of $\bP \ker(\pi_D)$ is the line spanned by
$\Delta(\cO^{\oplus 2})$ and $\Delta(\cO(p-q)\oplus\cO(q-p))$. This line
lies in $\bT_0$.

\item Suppose $h^0(\cO(2D)) = 1$ (this will be the case if $g \geq 4$
and $D \in C^{(2)}$ is general). Let $\epsilon$ be a point of
$\langle 2p+2q+r\rangle \setminus \bP Ker(\pi_D )$ and let
\[
0 \lra \cO(-D)\lra F_{\epsilon} \lra \cO(D) \lra 0
\]
be the extension defined by $\epsilon$. Then $h^0(F_{\epsilon}(D))=1$.

\item\label{crd} Choose coordinates $x,y,z,t$ on $\langle
2p+2q+r\rangle$ such that
\[
\bP \ker(\pi_D)=\;\{t=0\},\quad\langle 2p +
r\rangle =\;\{x=0\},\quad\langle 2q + r\rangle =\;\{y=0\},\quad\langle
p + q + r\rangle =\;\{z=0\}.
\]
Then the map $\psi_D$ restricted to $\langle 2p+2q+r\rangle$ factors
through a 4-dimensional linear subspace $\bP^4_{p+q,r}$ of $|2\Theta|$
and is given by:
\[
\psi_D:\langle 2p+2q+r\rangle\longrightarrow
\bP^4_{p+q,r}\;\hookrightarrow |2\Theta|,\qquad (x:y:z:t)\longmapsto
(xy:z^2:xt:yt:zt).
\]
\item\label{p4r}
The image $Y_{p+q,r}$ of $\langle 2p+2q+r\rangle$ by $\tpsi_D$ is a
cubic threefold in $\bP^4_{p+q,r}$ and for a good choice of coordinates:
\[
Y_{p+q,r} =\tpsi_D (\langle 2p+2q+r\rangle)=\quad\{ X_1X_2X_3-X_0X_4^2=0\}.
\]
The threefold $Y_{p+q, r}$ contains the points (we write $\Delta(D)$ for
$\Delta(\cO(D)\oplus\cO(-D))$):
\[
\begin{array}{ll}
\Delta(\cO)=(1:0:0:0:0),\qquad&
\Delta(p-q)=(0:1:0:0:0),\\
\Delta(p-r)=(0:0:1:0:0),\qquad&
\Delta(q-r)=(0:0:0:1:0).
\end{array}
\]
\end{enumerate}

\ts 1.$\quad$ From the definition of $\bP \ker(\pi_D)$ we see that its
defining (linear) equations are the elements of $H^0(\omega)\;(\subset
H^0(\omega(2D)))$. Since $H^0(\omega)\subset H^0(\omega(2p))\subset
H^0(\omega(2D))$ we see that $l_q\subset \bP \ker(\pi_D)$, similarly
$l_p\subset \bP \ker(\pi_D)$ and these two distinct lines already span
$\bP \ker(\pi_D)$. Similarly $H^0(\omega)\subset
H^0(\omega(p+q))\subset H^0(\omega(2D))$ verifies that the line $l$
lies in $\bP \ker(\pi_D)$.

\noindent
2.$\quad$ The map $\psi_D$ on $\bP \ker(\pi_D)$ is given by the
restriction of elements of $I_C(2)$ to the plane $\bP \ker(\pi_D)$. These
are conics passing through $p,\;q\;(\in \bP \ker(\pi_D))$, tangent to
$l_p$ at $p$ and tangent to $l_q$ at $q$. There is a pencil of conics
with these properties and since it contains $l^2$ and $l_pl_q$, it is
spanned by them.

\noindent
3.$\quad$ It follows from section \ref{psisym} that $\psi_D$ contracts
$l=<p,q>$ to the point $\Delta(\cO(p+q-D)\oplus
\cO(D-(p+q))=\Delta(\cO^2)$. Similarly:
$$
\psi_D(l_p)=\Delta(\cO(2p-D)\oplus \cO(D-2p))=
\Delta(\cO(p-q)\oplus \cO(q-p))=
\Delta(\cO(2q-D)\oplus \cO(D-2q)) =
\psi_D(l_q).
$$
By $2$, the image of $\bP Ker(\pi_D)$ is a line which therefore is
$\langle \Delta(\cO^{\oplus 2}),\Delta(\cO(p-q)\oplus \cO(q-p))
\rangle \subset \langle C-C \rangle = \bT_0$.

\noindent
4.$\quad$ Tensor the defining sequence of $F_\epsilon$ by $\cO(D)$, we
obtain
\[
0\longrightarrow \cO\longrightarrow F_\epsilon(D)\longrightarrow
\cO(2D)\longrightarrow 0.
\]
Now $h^0(F_\epsilon(D))\geq 2$ iff the coboundary map
$\delta_\epsilon:H^0(\cO(2D))\rightarrow H^1(\cO)$ is not injective. The
dual of the map
\[
H^1(\cO(-2D)) \lra \mbox{Hom}(H^0(\cO(2D)),H^1(\cO))\cong
H^0(\cO(2D))^*\otimes H^0(\omega)^*, \qquad\quad
\epsilon\longmapsto\delta_\epsilon
\]
is the multiplication map $m:H^0(\cO(2D))\otimes
H^0(\omega)\rightarrow H^0(\omega(2D))$. Since $h^0(\cO(2D))=1$, the
map $m$ is injective and its dual can be identified with
$\pi_D$. Therefore the map $\delta_{\epsilon}$ fails to be injective
exactly when $\epsilon$ is in $Ker(\pi_D)$.

\noindent
5.$\quad$ On $\langle 2p+2q+r\rangle$ the map $\psi_D$ is given by
quadrics passing through $p,\;q,\;r$ and whose restriction to
$\bP\ker(\pi_D)$ lies in $\langle xy,z^2\rangle$, thus they lie in the
space $\langle xy,\,z^2,\,xt,\,yt,\,zt\rangle$. For $\epsilon \in
\langle 2p+2q+r \rangle\setminus\langle 2p + 2q\rangle$, we have
$h^0(F_\epsilon(D))=1$ by 4. This implies that, up to multiplication
by a scalar, there is a unique nonzero map $\cO(-D)\rightarrow
F_\epsilon$ and gives the exact sequence defining the extension.
Hence the restriction of $\psi_D$ to the plane $\langle p+q+r \rangle$
is birational to its image. The conics in this plane passing through
$p,q,r$ span $\langle yt,\,xt,\,xy\rangle$. Since the map is
birational, it is given by these coordinate functions. Thus
$xy,\,z^2,\,xt,\,yt$ are coordinate functions. The only other possible
coordinate function is $zt$, if it is absent the map $\psi_D
|_{\langle 2p + 2q +r\rangle }$ is 2:1 (due to
the $z^2$), hence $zt$ is also a coordinate function.

\noindent
6.$\quad$ The coordinate functions obviously satisfy the cubic
equation, since the equation is irreducible and the image of $\tpsi_D$
is reduced (\cite{hartshorne77} page 92) and has dimension three it must
be the defining equation of the image.
\qed

\subsection{}
We now determine the relation between the linear subspaces $\bP^4_{p+q,r}$
and the hyperplanes $H_L$ and $H_W$ which we defined in the
introduction. We let $U_0$ be the open subset of $U'$ parametrizing
points $(p+q, r)$ such that the images of $p,q,r$ in $|\omega |^*$ are
not colinear.

\subsection{Proposition.} \label{prop411} Suppose $(p+q, r)\in U'$.
\begin{enumerate}
\item\label{pr1}
Suppose $L\in Sing(\Theta)$. Then: 
\[
\bP^4_{p+q,r}\subset H_L.
\]
\item \label{pr2}
Let $W$ be an element of $Gr(3,H^0(\omega))\setminus D_C$ and let
$H_W\;\in|2\Theta|^*$ be the corresponding hyperplane in $|2\Theta|$.
Then, for $(p+q, r)\in U_0$,
\[
\bP^4_{p+q,r}\subset H_W\;\Longleftrightarrow \; 
\langle\wedge^3W\cdot p\wedge q\wedge r\rangle=0.
\]
\end{enumerate}

\ts 1.$\quad$ Since $\bP^4_{p+q,r}$ is the span of $\psi_D(\langle
2p+2q+r\rangle)$, we must show that $F_\epsilon=\psi_D(\epsilon)$ lies
in $H_L$ for all $\epsilon \in \langle 2p+2q+r\rangle$.  So we must
show that $h^0(F_\epsilon\otimes L)\neq 0$. Tensoring the sequence
defining $F_\epsilon$ by $L$ we see that $L(-D)\hookrightarrow
F_\epsilon\otimes L$, so we are done if $h^0(L(-D))\neq 0$. Suppose
therefore that $h^0(L(-D))=0$, i.e., $h^0(\omega \otimes
L^{-1}(D))=2$. In particular, we have $h^0(L)=h^0(\omega\otimes
L^{-1})=2$. Since $H_L=H_{\omega\otimes L^{-1}}$ we may also assume
that $h^0(\omega \otimes L^{-1}(-D))=0$, i.e., $h^0(L(D))=2$.

It remains to show that the coboundary map
\[
\delta_\epsilon:H^0(L(D))\longrightarrow H^1(L(-D))\cong 
H^0(\omega\otimes L^{-1}(D))^*,
\]
with $L\in Sing(\Theta),\;h^0(L(D))=h^0(\omega\otimes L^{-1}(D))=2$,
is not injective. 

The dual of the map $\epsilon\mapsto \delta_\epsilon$ is again the 
multiplication map:
\[
m:H^0(\omega\otimes L^{-1}(D))\otimes H^0(L(D))\longrightarrow
H^0(\omega(2D))= H^1(\cO(-2D))^*.
\]
Choose bases $\{ s_1',\,s_2'\}$ of $H^0(L)$ and $\{ t_1',\,t_2'\}$ of
$H^0(\omega\otimes L^{-1})$ and let $\{ s\}$ be a basis of
$H^0(\cO(D))$.  Let $h_{ij}:=m(ss_i'\otimes st_j')\;\in
H^1(\cO(-2D))^*$. Then $\delta_\epsilon$ is not injective iff
$Q(\epsilon)=0$ with $Q:=h_{11}h_{22}-h_{12}h_{21}$. Note that $Q =
s^4 \oQ$ where $\oQ := g_{11}g_{22}-g_{12}g_{21}$ with $g_{ij} :=
h_{ij}/s^2$. So $\oQ$ is the equation of the tangent cone to $\Theta$
at the singular point $L$ (\cite{ACGH} page 240). Therefore $\oQ$
defines a quadric $\oq$ in $|\omega|^*$ which contains $C$
(\cite{ACGH} page 241) and the inverse image of $\oq$ in
$|\omega(2D)|^*$ is the quadric $q$ of equation $Q$ which contains
$C_D$ and whose vertex contains $\bP Ker(\pi_D) = \langle 2p+2q
\rangle$. It follows that $q$ contains $\langle 2p+2q+r\rangle$ for
all $r \in C$ and hence $\bP^4_{p+q,r} \subset H_L$.

\noindent
2.$\quad$ The assumption on the points implies that $(p\wedge q\wedge
r)^\perp= H^0(\omega(-p-q-r))\;\subset H^0(\omega)$ has dimension
$g-3$ and $\langle\wedge^3W\cdot p\wedge q\wedge r\rangle=0$ is
equivalent to $W\cap H^0(\omega(-p-q-r))\neq \{0\}$.  We will write
$D=p+q$. A result of Beauville (\cite{beauville91} page 268) on the
intersection of $H_E:=\D_\omega^*(E)$ and $\cM_{\cO}$ is:
\[
\Delta(F)\in H_E\Longleftrightarrow h^0(E\otimes F)\neq 0
\qquad\qquad F\in \cM_\cO,\;E\in\cM_\omega.
\]
Thus we have to show:
\[
(h^0(E_W\otimes F_\epsilon)\neq 0
\quad\forall \epsilon\in \langle 2p+2q+r\rangle)\;
\Longleftrightarrow\; (W\cap H^0(\omega(-p-q-r))\neq \{0\}).
\]
We will write $E$ for $E_W$.
The defining sequence of $F_\epsilon$ tensored by $E$ shows
$H^0(E(-D))\hookrightarrow H^0(E\otimes F_{\epsilon})$. If $h^0(E(-D))
> 0$, then $h^0(E\otimes F_\epsilon)\neq 0$ for all $\epsilon\in
\langle 2p+2q+r\rangle$. On the other hand
$h^0(E(-D))=h^0(\omega\otimes E^*(-D))=h^1(E(D))=h^0(E(D))-4$ by
Riemann-Roch, Serre duality and the fact that
$E\cong\omega\otimes E^*$. Thus $h^0(E(-D)) > 0
\Leftrightarrow h^0(E(D)) > 4$. By Lemma \ref{secEt},
$h^0(E(D))\geq 5$ implies $\dim W_E\cap H^0(\omega(-p-q))>1$ and thus $\dim
W_E\cap H^0(\omega(-p-q-r))>0$.

It remains to consider the case where $h^0(E(-D)) = 0$ and $h^0(E(D))=4$.
Now $h^0(E\otimes F_\epsilon)\neq 0$ iff the coboundary map
\[
\delta_\epsilon:H^0(E(D))\longrightarrow H^1(E(-D))\cong H^0(E(D))^*
\]
is not injective.

Let $\{ t_1,\ldots,\,t_4\}$ be a basis of $H^0(E(D))$. Let
$\alpha:\wedge^2H^0(E(D))\rightarrow H^0(\wedge^2
E(D))=H^0(\omega(2D))\cong H^1(\cO(-2D))^*$ be the evaluation map and
define $h_{ij}:=\alpha(t_i\wedge t_j)\;(\in H^1(\cO(-2D))^*)$. Let
$P_E$ be the pfaffian (a square root of the determinant) of the
$4\times 4$ alternating matrix $(h_{ij})$. Thus $P_E$ is an element of
$S^2H^0(\omega(2D))$ and defines a quadratic form on $H^1(\cO(-2D))$
also denoted by $P_E$. By construction $\Delta(F_\epsilon)\in H_E$ iff
$P_E(\epsilon)=0$. Thus we must show:
\[
{P_E}|_{\langle 2p+2q+r\rangle}=0\;\Longleftrightarrow\; W\cap
H^0(\omega(-p-q-r))\neq \{0\}.
\]
We now choose the basis $\{ t_i\}$ of $H^0(E(D))$ with some more
care. Let $s_1\in H^0(E)$ be a generator of the (one-dimensional)
kernel of the evaluation map $H^0(E)\rightarrow E_p$. Then $s_1(q)\neq
0$ because $h^0(E(-D)) = 0$. Let $s_2\in H^0(E)$ be a generator of
$Ker(H^0(E)\rightarrow E_q)$ and let $s_3$ be a third section such
that $\{ s_1, s_2, s_3 \}$ is a basis of $H^0(E)$. Then $s_3$ is
not zero at $p$ nor $q$. Put
\[
t_i:=ss_i\quad\in H^0(E(D)),
\]
where $\{ s\}$ is a basis of $H^0(D)$.


Next let $W=\{t\in H^0(E(D)):\,t(r)=0\}$, a vector space of dimension 
at least $2$. Since $r\not\in supp(D)$ and $E$ is generated by 
global sections, the image of 
$H^0(E)\stackrel{s}{\hookrightarrow} H^0(E(D))$ intersected with $W$ 
has dimension at most one. Hence there exits a $t_4\in H^0(E(D))$ with
$t_4(r)=0$ and $t_4\neq st$ for any $t\in H^0(E)$.

Now we observe that $\alpha(t_1\wedge t_4)=\alpha(ss_1\wedge t_4)$ is
zero in $p$ with multiplicity $\geq 2$ (since $s(p)=0$ and $s_1(p)=0$)
and is also zero in $q$ ($s(q)=0$) and $r$ ($t_4(r)=0$). Therefore
this element of $H^0(\omega(2D))$ must be zero on the subspace
$\langle 2p+2q+r\rangle$ (which is spanned by $q,\,r$ and $l_p
=\langle 2p\rangle$ the tangent
line to $C_D$ at $p$). Similarly $\alpha(t_2\wedge t_4)$ must be zero
at $p+2q+r$ and thus it is zero on the subspace $\langle
2p+2q+r\rangle$. The pfaffian $P_E$,
restricted to $\langle 2p+2q+r\rangle$ is then easy to compute:
\[
{P_E}|_{\langle 2p+2q+r\rangle} =\alpha(t_1\wedge t_2 )\odot\alpha(
t_3\wedge t_4 )\qquad \in S^2 H^0 (\omega( 2D )).
\]

We claim:
\[
\alpha(t_1\wedge t_2)|_{\langle 2p+2q+r\rangle}=0\;\Longleftrightarrow\;
\dim W_E\cap H^0(\omega(-p-q-r))=1
\]
and 
\[
\alpha(t_1\wedge t_2)|_{\langle 2p+2q+r\rangle} \neq
0\;\Longrightarrow\;
\alpha(t_3\wedge t_4)|_{\langle 2p+2q+r\rangle}\neq 0.
\]
This will prove 2.

Since $h^0(E(D))=4$ we have $\dim W_E \cap H^0(\omega(-p-q))=1$ by
Lemma \ref{secEt}. By construction $0\neq\alpha(s_1\wedge s_2)\;\in
H^0(\omega)$ is zero in $p$ and $q$. Thus $\langle \alpha(s_1\wedge
s_2)\rangle =W_E\cap H^0(\omega(-p-q))$. This implies:
\[
\dim W_E\cap H^0(\omega(-p-q-r))=1 \Longleftrightarrow
\alpha(s_1\wedge s_2)(r)=0.
\]
By construction, $\alpha(t_1\wedge t_2)$ has double zeros in $p$ and
$q$; if it also has a zero in $r$ then it vanishes on $\langle
2p+2q+r\rangle$. Conversely, if $\alpha(t_1\wedge t_2)|_{\langle
2p+2q+r\rangle}=0$, then $\alpha(t_1\wedge t_2)=\alpha(ss_1\wedge
ss_2)$ must vanish in $r$. Since $s(r)\neq 0$, we have $\alpha(s
s_1\wedge s s_2)(r)= 0\Leftrightarrow\alpha( s_1\wedge s_2)(r)=
0$.

It remains to show that the restriction of $\alpha(t_3\wedge t_4)$ to
$\langle 2p+2q+r\rangle$ is non-zero if $\alpha(t_1\wedge
t_2)|_{\langle 2p+2q+r\rangle} \neq 0$. So assume $\alpha(t_1\wedge
t_2)|_{\langle 2p+2q+r\rangle} \neq 0$. Then $s_1\wedge s_2$ is not
zero at $r$ and we can take $s_3$ to be a basis of $Ker(H^0(E) \lra
E_r)$. The differential $\alpha(t_3\wedge t_4)$ is zero in $p$ and $q$
(since $t_3=ss_3$ and $s$ is zero there) and also in $r$, thus we must
show it does not have a double zero at $q$ (or, equivalently, at $p$).

Let $F_{x}$ be the stalk of a sheaf $F$ on $C$ at $x\in C$. From the choice 
of the $s_i\in H^0(E)$, we can define an isomorphism $\phi:E_{q}\cong 
\cO^{\oplus 2}_{q}$ such that, with $z$ a local parameter at $q$:
\[
\phi:\quad
s_1\longmapsto (1,0)+z(\ldots),\qquad s_2\longmapsto z(a,b)+z^2(\ldots),
\qquad
s_3\longmapsto (0,1)+z(\ldots).
\]
Also, we can define an isomorphism $\psi:E(D)_{q}\cong\cO^{\oplus
2}_{q}$, such that $\psi(t_i)=z\phi(s_i)$ for $i \in\{ 1,2,3\}$. The
section $t_4$ is not zero at $q$, thus we can write
\[
\psi(t_4)=(c,d)+z(\ldots)\;\in\cO^{\oplus 2}_{x}\qquad{\rm  with}\quad 
(c,d)\neq (0,0)\quad\in\bC^2.
\]
Now we recall that $\alpha(t_1\wedge t_4)$ has a double zero in $q$.
In $q$ the local expansion of $t_1\wedge t_4$ is:
\[
(t_1\wedge 
t_4)_q=\left(z(1,0)+z^2(\ldots)\right)\wedge
\left((c,d)+z(\ldots)\right)=dz+z^2(\ldots)
\]
thus we must have $d=0$ and we see that $c\neq 0$. Therefore:
$$
(t_3\wedge t_4)_q=\left(z(0,1)+z^2\ldots)\right)\wedge
\left((c,0)+z(\ldots)\right)=
cz+z^2(\ldots).
$$
We conclude that $\alpha(t_3\wedge t_4)$ has a zero of order exactly one 
at $q$.
This proves the claim and completes the proof of the proposition. \qed

\subsection{Corollary.}\label{capT0}
For any $(p+q, r)\in U'$ and $W\in Gr(3,H^0(\omega))\setminus D_C$ we have:
\[
\bP^4_{p+q,r}\subset \bT
\]
and
\[
\bT_0\subset H_W.
\]
If, moreover, $(p+q,r)\in U_0$,
\[
\dim \bP^4_{p+q,r}\cap \bT_0=3.
\]
For $W\in Gr(3,H^0(\omega)) \setminus D_C$,
\[
\bT\not\subset H_W.
\]

\ts The intersection of $\bP^4_{p+q,r}$ with $\Delta(\cM_{\cO})$ contains
the cubic threefold determined in \ref{theta}. Since this threefold is
singular at $\Delta(\cO^{\oplus 2})$, its embedded tangent space at
$\Delta(\cO^{\oplus 2})$ is all of $\bP^4_{p+q,r}$. Thus $\bP^4_{p+q,r}$ is
contained in $\bT$.

Since $h^0(E_W)= 3$, $h^0(E(-p))\geq 1$ for all $p\in C$. Therefore
$C-C\;\subset Pic^0C$ is contained in $\Delta_\omega(E_W)$ and so
its span, which is $\bT_0$ (cf.\ \ref{C-Cperp}), is contained in $H_E=H_W$.

Since $\Delta(\cO),\;\Delta(p-q),\;\Delta(p-r),\;\Delta(q-r)$ are in
$\bP^4_{p+q,r}$ (see \ref{theta}) and in the span of $C-C$, they lie in
$\bT_0\cap \bP^4_{p+q,r}$. As they are independent (cf.\ \ref{theta})
we obtain $\dim
\bT_0\cap \bP^4_{p+q,r}\geq 3$. For $p,q,r \in C_{can}$ distinct and
non-colinear, the subvariety $\langle p+q+r \rangle^{\perp}$ of
$Gr(3,H^0(\omega))$ consisting of those $W \in Gr(3,H^0(\omega))$ such
that $\langle \wedge^3W,p\wedge q\wedge r\rangle = 0$ is a divisor in
the Pl\"ucker system. Choosing $W \in Gr(3,H^0(\omega)) \setminus
\left( D_C \cup \langle p+q+r \rangle^{\perp}\right)$, we have $\langle
\wedge^3W,p\wedge q\wedge r\rangle\neq 0$. Then
$\bP^4_{p+q,r}\not\subset H_W$ by Proposition \ref{prop411}. However,
$\bT_0\subset H_W$ hence $\bP^4_{p+q,r}\not\subset\bT_0$.

As $p,q$ and $r$ vary, the divisors $\langle p+q+r \rangle^{\perp}$
span the Pl\"ucker linear system on $Gr(3, H^0 (\omega ))$ which has
no base points. Hence, given $W\in Gr(3,H^0(\omega))\setminus D_C$,
we can find $p,q,r$ with $\langle \wedge^3W,p\wedge q\wedge
r\rangle\neq 0$. Then, by Proposition \ref{prop411}, we have
$\bP^4_{p+q,r}\not\subset H_W$. Since $\bP^4_{p+q,r}\subset \bT$ we
obtain $\bT\not\subset H_W$. \qed
\vskip10pt

Let $U''$ be the open subset of $C^{ (3)}$ parametrizing divisors $E$
such that $h^0 (E) = 1$. Note that, since $C$ is not hyperelliptic,
either $U'' = C^{ (3)}$ or the complement of $U''$ has codimension two
in $C^{ (3)}$. We have

\subsection{Proposition.}\label{proppl}
The morphism
\[
\rho: U_0\lra \bT/\bT_0,\qquad
(p+q,r)\longmapsto \bP^4_{p+q,r}/(\bP^4_{p+q,r}\cap\bT_0)
\]
is the composition of the morphism $\pi: U_0\inj C^{(2)}\times
C\rightarrow C^{(3)},\quad (p+q,r)\mapsto p+q+r$ with the morphism:
\[
\begin{array}{cccccc}
P:& U'' &\lra &Gr(3,H^1(\cO)) &\longrightarrow& \bP(\wedge^3H^1(\cO)),\\
& p+q+r&\longmapsto&\langle p+q+r\rangle& \longmapsto
&\wedge^3\langle p+q+r\rangle,
\end{array}
\]
and with a linear embedding $\bP(\wedge^3H^1(\cO))\hookrightarrow \bT/\bT_0$.

In particular we have:
\[
\dim \langle \cup \bP^4_{p+q,r}\rangle =\half{g(g+1)} + {g\choose 3}.
\]

Furthermore, the morphism
\[
Gr(3,H^0(\omega))\setminus D_C\longrightarrow
(\bT/\bT_0)^*,\qquad W\longmapsto \bar{H}_W
\]
is the composition of the Pl\"ucker map with a linear embedding
$\bP(\wedge^3H^0 (\omega ))\hookrightarrow (\bT/\bT_0)^*$.

\ts 
By Corollary
\ref{capT0} any $W\in Gr(3,H^0(\omega))\setminus D_C$ gives a hyperplane
$\bar{H}_W:=(H_W\cap\bT)/\bT_0$ in $\bT/\bT_0$. We first want to determine
the divisor $D_W:=\rho^*\bar{H}_W$. 
From \ref{prop411}.\ref{pr2} we
have the following equality of sets:
\[
D_W=\{(p+q,r):\;\langle\wedge^3W,p\wedge q\wedge r\rangle=0\;\}.
\]
On the other hand, let $D'_W$ be the inverse image under $P$ of the
hyperplane 
$(\wedge^3W)^\perp\;\subset\bP(\wedge^3H^1(\cO))$, i.e.,
$D'_W:=P^*(\wedge^3W)^\perp$.  Then, again as sets,
\[
D'_W=P^{-1}((\wedge^3W)^\perp)=\{p+q+r:\;\langle\wedge^3W,p\wedge
q\wedge r\rangle=0\;\}.
\]
Hence $D_W =\pi^{-1}D'_W$ as sets. 

In Lemma \ref{redirr} below we prove that the divisors $D'_W$ and
their inverse images in $C^3$ are reduced and irreducible for general
$W$. Hence $\pi^*(D'_W)\subset U_0$ is reduced and irreducible. Thus
$D_W = n\pi^*D'_W$ as divisors for some positive integer $n$. The
divisor $D_W'$ is the zero locus of $s_W\in
H^0(U'',P^*\cO(1)) = H^0( C^{(3) },P^*\cO(1))$. The map
$\pi^*:H^0( U'' ,P^*\cO(1))\rightarrow H^0(U_0,\rho^*\cO(1))$ is
linear on the one hand and maps $s_W\mapsto (\pi^*(s_W))^{\otimes n}$
on the other hand. We conclude that $n=1$ and that $\rho$ is given by
sections of the Pl\"ucker line bundle.

Now we need to prove that $\rho$ is given by the full Pl\"ucker linear
system. For this, we need to produce ${g\choose 3}$ elments of $U_0$
whose images by $\rho$ are linearly independent. Since $P(\pi( U_0
))$ spans $\bP(\wedge^3H^1(\cO))$ (just take $g$ general
points on $C_{can}\subset \bP H^1(\cO)$ then the wedges of three
distinct points are independent), we can find $D_i\in \pi( U_0)$ and
$W_i\in Gr(3,H^0(\omega))\setminus D_C$ (for $i=1,\ldots,{g\choose
3}$) such that $P(D_i)\not\in (\wedge^3W_i)^\perp$ and $P(D_i)\in
(\wedge^3W_j)^\perp$ if $j\neq i$. For any choice of $(p_i+q_i,r_i)\in
\pi^{-1}(D_i)$, the images $\rho(p_i+q_i,r_i)$ are independent.

The last statement follows from the natural duality between $G(3,
H^1(\cO ))$ and $G(3, H^0(\omega ))$ and the fact that $\rho$ is
the composition of $\pi$ with the Pl\"ucker map.\qed

\subsection{Lemma.} \label{redirr}
Let $W$ be a general three-dimensional subvector space of
$H^0(\omega)$. Then the divisor $D'_W :=\{ p+q+r :\;\langle\wedge^3
W,p\wedge q\wedge r\rangle=0\;\}$ in $U''$ and its inverse image
in $C^3$ are irreducible and reduced.

\ts $\quad$ Set-theoretically, the divisor $D'_W$ parametrizes the
divisors $D\in U''$ such that there exists a divisor $D'\in
|W|\subset|\omega|$ with $D \leq D'$ (i.e., the divisor $D'-D$ is
effective). The cohomology class of the closure of this subset in
$C^{(3)}$ is the coefficient of $t$ in the expression
(\cite{macdonald62}, 16.2, p.\ 338)
\[
(1+\eta.t)^{g-4}\prod_{i=1}^{g}(1+A_{i}B_{i}.t),
\]
i.e.,
\begin{equation}
(g-4)\eta+ \sum_{i=1}^{g}A_{i}B_{i}.\label{class}
\end{equation}
To determine the cohomology class of $D'_W$ we proceed as follows.
Let $D\subset C^{(3)}\times C$ be the universal divisor, then we have
an exact sequence:
\[
0\lra \cO_{C^{(3)}\times C}\lra \cO_{C^{(3)}\times C}(D)\lra
\cO_D(D)\lra 0.
\]
Let $p: C^{(3)}\times C\rightarrow C^{(3)}$ be the first projection and
consider $p_*$ of the sequence above. It is easy to see that $p_*\cO_{
C^{(3)}\times C}\cong p_*\cO_{ C^{(3)}\times C} (D)\cong\cO_{
C^{(3)}}$ and, by the general theory of Hilbert schemes
(\cite{grothendieck60}), we have $p_*\cO_D(D)\cong T_{C^{(3)}}$.
Furthermore, since $p |_D$ is finite, all higher direct images by $p$
of $\cO_D (D)$ are zero. So we have the short exact
sequence of sheaves 
(use $p_*\cO_D (D)\cong T_{C^{(3)}}$ \cite{grothendieck60}):
\[
0\lra T_{C^{(3)}}\lra H^1(\cO)\otimes\cO \lra Q\lra 0
\] 
where we have put $Q := R^1 p_*\cO_{ C^{(3)}\times C} (D)$. For every
$E\in U''$, the image of the fiber $(p_*\cO_D(D))_E$ in $\bP H^1(\cO)$
is $\langle E\rangle$, hence $T_{C^{(3)}} |_{U''}=\rho^*{\bf W}$ where
$\bf W$ is the universal subbundle on $Gr(3,H^1(\cO))$. Since the
Pl\"ucker map is given by $\wedge^3 {\bf W}^*$, the cohomology class
of $D'_W$ is $-c_1(\wedge^3 T_{C^{(3)}})=c_1(\omega_{C^{(3)}})$. The
class $c_1(\omega_{C^{(3)}})$ is determined in \cite{macdonald62},
(14.9), p.\ 334 and \cite{ACGH}, VII, (5.4) and is equal to the class
(\ref{class}) above.

Since the two cohomology classes are equal, we deduce that the two
divisors are equal as schemes and $D'_W$ is reduced. Hence its inverse
image in $C^3$ is also reduced. The irreducibility of the inverse
image of $D'_W$ in $C^3$ (and hence the irreducibility of $D'_W$)
follows from the Lemma on page 111 of \cite{ACGH} (the argument there
works for $r=2$ as well, there is a misprint in the statement of the
Uniform Position Theorem on page 112).
\qed

\subsection{Proof of Theorem \ref{emb}}\label{dDinj} 
It remains (cf.\ Corollary \ref{prfgenxi}) to consider the case
$\xi^{\otimes 2}\cong\cO$. Since the points of order two of $Pic^0C$
act on $\cM_\cO$ (by $F\mapsto F\otimes\xi$)
and on $|2\Theta|$ ($D\mapsto D+\xi$)
it suffices to consider the case $\xi=\cO$.
Thus we need to show that the differential
$({\rm d}\Delta)_{\cO^{\oplus 2}}$ of
$\Delta:\cM_\cO\rightarrow |2\Theta|$ at $\cO^{\oplus 2}$ is injective. 

Choose $( p+q, r)\in U_0$. Recall that $Y_{p+q,r}\subset\bP^4_{p+q,r}$
is the image of $\langle 2p+2q+r\rangle\cong\bP^3$ by $\tpsi_{p+q}
:=\D\tphi_{p+q}$. By Proposition \ref{theta}.6, there is a choice of
coordinates $X_0,\ldots,X_4$ on $\bP^4_{p+q,r}$ such that $Y_{p+q,r}$
has equation $X_1X_2X_3-X_0X_4^2=0$. A simple computation shows that
$Y_{p+q,r}$ is nonsingular in codimension $1$. Since $Y_{p+q,r}$ is
Cohen-Macaulay (complete intersection in a smooth scheme), it follows
that $Y_{p+q,r}$ is normal. Put $X_{p+q,r}:=\tphi (\langle 2p + 2q
+r\rangle )\;\subset\cM_\cO$. Then $X_{p+q,r}$ is reduced
(\cite{hartshorne77} page 92). The morphism $\D$ is injective
(\cite{brivioverra96}), thus $\D|_{X_{p+q,r}}:X_{p+q,r}\rightarrow
Y_{p+q,r}$ is bijective. As $Y_{p+q,r}$ is normal, the morphism
$\D|_{X_{p+q,r}}$ is an isomorphism.

Thus $({\rm d}\D)_{\cO^{\oplus 2}}$ induces an isomorphism between
$T_{p+q,r} := T_{\cO^{\oplus 2}} X_{p+q,r}$ (a subspace of
$T_{\cO^{\oplus 2}}\cM_\cO$) and $T_{\D(\cO^{\oplus 2})}Y_{p+q,r}$, a
4-dimensional vector space. Hence $({\rm d}\D)_{\cO^{\oplus 2}}
({T_{p+q,r}}) = T_{\Delta(\cO^{\oplus 2})}\bP^4_{p+q,r}$. Proposition
\ref{proppl} shows $\dim\langle\cup\bP^4_{p+q,r}\rangle= \dim
T_{\cO^{2}}\cM_{\cO}$. Therefore the dimension of the image of $({\rm
d}\D)_{\cO^{\oplus 2}}$ is at least $\dim T_{\cO^{\oplus 2}}\cM_{\cO}$
which implies the injectivity of $({\rm d}\D)_{\cO^{\oplus 2}}$. \qed

\subsection{Corollary.}\label{TsubTheta}
We have
\[
\dim\bT=\dim T_{\cO^{\oplus 2}}\cM_\cO=\half{g(g+1)} + {g\choose 3}
\qquad{\rm and}\quad
\langle\cup\bP^4_{p+q,r}\rangle =\bT\;\subset\;\langle Sing(\T)\rangle^\perp.
\]

\ts
The first equalities follow from the fact that $\Delta$ is an embedding.
Next, by Theorem \ref{emb}, we have that 
\[
\dim \bT=\dim T_{\cO^{\oplus 2}}\cM_\cO=\half g(g+1)+
{g\choose 3}.
\]
This is the same as the dimension of $\langle\cup\bP^4_{ p+q,r }\rangle$
by Proposition \ref{proppl}. Since each $\bP^4_{p+q,r}\subset \bT$
(cf.\ \ref{capT0}), we obtain $\langle\cup\bP^4_{p+q,r}\rangle =\bT$. 
The inclusion $\bT\subset\;\langle Sing(\T)\rangle^\perp$
follows
from the fact that each $\bP^4_{p+q,r}\subset H_L$ for any $L\in
Sing(\Theta)$ (cf.\ \ref{prop411}.\ref{pr1}), hence $\bT\subset \cap H_L=\langle
Sing(\T)\rangle^\perp$.  \qed

\subsection{Proof of Theorem \ref{t2}}\label{prt2}
By Proposition \ref{proppl}, the map:
\[
Gr(3,H^0(\omega))\setminus D_C\longrightarrow
(\bT/\bT_0)^*,\qquad W\longmapsto \bar{H}_W
\]
where $\bar{H}_W\;\subset\bT/\bT_0$ is the image of $H_W\cap
\bT$, is the Pl\"ucker map. For dimension reasons we then get 
the isomorphism $\bP (\wedge^3 H^0(\omega) )\cong (\bT /\bT_0)^*$.
\qed

\section{The span of $Sing(\Theta_{\xi})$}\label{span}

\subsection{}
In this section we prove Theorem \ref{mainthm}.
We have already seen in \ref{TsubTheta} and \ref{prfgenxi} that
\[
\bT_{\xi}\;\subset\;\langle Sing(\T_{\xi})\rangle^\perp,
\]
and that the dimensions of the embedded tangent spaces are:
\[
\dim\bT_{\xi }=\half {g(g+1)} + {g\choose 3}\;\hbox{  if }\xi^{\otimes
2}\cong\cO ,\qquad\dim\bT_{\xi }= g^2 - g + 1\;\hbox{ if }\xi^{\otimes
2}\not\cong\cO .
\]
Thus to prove Theorem \ref{mainthm} it suffices to prove 
that the dimension of the linear subspace
$\langle Sing(\Theta_{\xi})\rangle^{\perp}$ of $|2\Theta|$ is
equal to the dimension of $\bT_{\xi}$. This subspace
is the linear system of $2\T$-divisors
which contain $Sing(\T_{\xi})$. By translation by $\xi$, it can be
identified with the linear system of $2\T_{\xi}$-divisors which
contain $Sing (\T)$, i.e., the projectivization of
$H^0(\cI\otimes\cO_{Pic^{g-1}C}(2\T_{\xi}))$ where $\cI$ is the ideal sheaf
of $Sing(\T)$ in $Pic^{g-1}C$. Put
\[
\langle Sing(\T)\rangle^\perp_{\xi} :=
\{D\in |2\T_{\xi}|:\;Sing(\T)\subset D\}=
\bP H^0(Pic^{g-1}C,\cI\otimes\cO_{Pic^{g-1}C}(2\T_{\xi})).
\]
We have the usual exact sequence
\[
0\lra\cO_{Pic^{g-1}C} ( 2\T_{\xi} -\T )\lra\cO_{Pic^{g-1 }C}
(2\T_{\xi })\lra\cO_{\T} (2\T_{\xi})\lra 0\: .
\]
By the theorem of the square,
\[
\cO_{Pic^{g-1}C} ( 2\T_{\xi} )\cong\cO_{Pic^{g-1}C} ( \T_{\xi^{\otimes
2}} +\T ) .
\]
Hence the above exact sequence becomes
\[
0\lra\cO_{Pic^{g-1}C} ( \T_{\xi^{\otimes 2}} )\lra\cO_{Pic^{g-1 }C}
(2\T_{\xi })\lra\cO_{\T} (2\T_{\xi})\lra 0\: .
\]
As $h^1(\cO_{Pic^{g-1}C} (\T_{\xi^{\otimes 2}} ))=0$, we have
$h^0(\cO_{\T} (2\T_{\xi })) = h^0(\cO_{Pic^{g-1}C} (2\T_{\xi
}))-1$.

Since the restriction of $\cI$ to $\T$ is the ideal sheaf of
$Sing(\T)$ in $\T$, we need to prove that
$h^0(\T,\cI\otimes\cO_{\T}(2\T_{\xi }) ) =\dim T_{\xi\oplus\xi^{
-1}}\cM_\cO$.

\subsection{}
To determine $H^0(\Theta,\cI\otimes\cO_\Theta(2\Theta_{\xi }))$ we use the
resolution of singularities of $\Theta$ obtained from the natural
morphism $C^{(g-1)}\rightarrow Pic^{g-1} C$:
\[
\rho : C^{(g-1)} \lra \T,\qquad Z:=\rho^{-1}(Sing(\Theta))\quad\subset
C^{(g-1)},\qquad L_{\xi } :=\rho^*\cO_{\T}(\T_{\xi }),\qquad L := L_{\cO }.
\]
Let $\cI_Z$ be the ideal sheaf of $Z$ in $C^{(g-1)}$. Then $\cI_Z
=\rho^*\cI$ (see \cite{ACGH} Proposition 3.4 page 181). We have

\subsection{Lemma.} The map $\rho$ induces an isomorphism:
\[
\rho^*:H^0(\Theta,\cI\otimes\cO_\Theta(2\Theta_{\xi }))
\lra H^0(C^{(g-1)}, \cI_Z \otimes L_{\xi }^{\otimes 2}).
\]

\ts Clearly $\rho^*$ is injective. The two spaces have the same
dimension since $H^0( C^{ (g-1) },\cI_Z \otimes L_{\xi }^{\otimes 2}) \cong
H^0(\T ,\rho_*(\cI_Z \otimes L_{\xi }^{\otimes 2}))\cong
H^0(\T, \cI\otimes\cO_{\T} (2\T_{\xi })\otimes\rho_*\cO_{C^{(g-1)}})$ (the
last isomorphism is obtained from the projection formula) and
$\rho_*\cO_{C^{(g-1)}} =\cO_{\T}$ because the fibers of $\rho$ are
connected. \qed

\subsection{} Therefore we need to determine the dimension of
$H^0(C^{(g-1)}, \cI_Z\otimes L_{\xi }^{\otimes 2})$. The ideal sheaf $\cI_Z$
has the resolution (see \cite{green84} page 88, \cite{ACGH}, VI.4, p.\ 258)
\[
0 \lra T_{C^{(g-1)}} \otimes L^{-1} \lra H^1(C, \cO)
\otimes L^{-1} \lra \cI_Z \lra 0.
\]
Tensoring the sequence with $L_{\xi }^{\otimes 2}$ and using the
theorem of the square, we obtain
\begin{equation}
0 \lra T_{C^{(g-1)}} \otimes L_{\xi^{\otimes 2} } \lra H^1(C, \cO_{ C
})\otimes L_{\xi^{\otimes 2} } \lra \cI_Z \otimes L_{\xi }^{\otimes 2}
\lra 0 \label{resolution}
\end{equation}

\subsection{Theorem.}\label{thmsing}
\begin{enumerate}
\item Suppose that $\xi^{\otimes 2}\not\cong\cO$. With the above
notation we have the exact sequence
\[
0 \lra H^1(C, \cO) \lra H^0(C^{(g-1)}, \cI_Z \otimes L_{\xi
}^{\otimes 2}) \lra H^1(C^{(g-1)}, T_{C^{(g-1)}} \otimes
L_{\xi^{\otimes 2} }) \lra 0
\]
and $h^1(C^{(g-1)}, T_{C^{(g-1)}} \otimes L_{\xi^{\otimes 2} }) =
(g-1)^2$. Thus 
\[
\dim\langle Sing(\T_{\xi })\rangle^{\perp} = g + (g-1)^2 = g^2 - g + 1,
\]
\[
\dim \langle Sing(\Theta_{\xi })\rangle=
2^g- 2 - g^2 + g
\]
and $\bT_{\xi }=\langle Sing(\T_{\xi })\rangle^\perp$.

\item Suppose that $\xi^{\otimes 2}\cong\cO$. We have the exact sequence
\[
0 \lra S^2H^1(C, \cO) \lra H^0(C^{(g-1)}, \cI_Z \otimes L_{\xi
}^{\otimes 2}) \lra\Lambda^{3}H^1(C, \cO) \lra 0.
\]
Thus
\[
\dim\langle Sing(\T_{\xi })\rangle^{\perp} = {{g+1}\choose 2}+{g\choose 3} = 
\sum_{i=1}^3 {g\choose i},
\]
\[
\dim \langle Sing(\Theta_{\xi })\rangle=
2^g-2-\left({{g+1}\choose 2}+{g\choose 3}\right)=
2^g-1-\sum_{i=0}^3 {g\choose i}
\]
and  $\bT_{\xi }=\langle Sing(\T_{\xi })\rangle^\perp$.

\end{enumerate}

\subsection{}
The exact sequences in the theorem are obtained from the cohomology sequence
of sequence (\ref{resolution}). Therefore, to prove the theorem,
 we will compute the cohomology vector spaces of
the sheaves appearing in sequence (\ref{resolution}).

\subsection{Lemma.}\label{Hi}
\begin{enumerate}
\item Suppose that $\xi^{\otimes 2}\not\cong\cO$. Then, for all $i >
0$,
\[
H^i(C^{(g-1)}, L_{\xi^{\otimes 2} }) = 0
\]
and $H^0(C^{(g-1)}, L_{\xi^{\otimes 2} })\cong H^0 (Pic^{g-1} C,\cO_{
Pic^{g-1} C} (\T_{\xi^{\otimes 2 } }))\cong\bC$.

\item For all $i\geq 0$
\[
H^i(C^{(g-1)}, L) \cong H^i(\T,
\cO_{\T}(\T)) \cong H^{i+1}(Pic^{g-1}C, \cO_{Pic^{g-1}C}) \cong
\Lambda^{i+1} H^1(C, \cO).
\]
\end{enumerate}

\ts Since the morphism $\rho : C^{(g-1)}\rightarrow\T$ is a rational
resolution, we have
\[
H^i(C^{(g-1)}, L_{\xi^{\otimes 2} })\cong H^i(\T,\cO_{\T
}(\T_{\xi^{\otimes 2 } }))
\]
for all $i\geq 0$ and all $\xi$. Now the lemma follows from the cohomology sequence of the usual exact
sequence
\[
0\lra\cO_{Pic^{g-1}C} (\T_{\xi^{\otimes 2} } -\T )\lra\cO_{Pic^{g-1}C}
(\T_{\xi^{\otimes 2} })\lra\cO_{\T} (\T_{\xi^{\otimes 2} })\lra 0\: .
\]
\qed
\vskip20pt

The above Lemma determines the cohomology of $H^1 (C,\cO)\otimes
L_{\xi^{\otimes 2} }$. Our next step is

\subsection{Lemma.}\label{H0TL}
We have
\[
H^0(C^{(g-1)},T_{C^{(g-1)}} \otimes L)\cong \Lambda^2H^1(C,\cO)
\]
and, if $\xi^{\otimes 2}\not\cong\cO$, then
\[
H^0(C^{(g-1)},T_{C^{(g-1)}} \otimes L_{\xi^{\otimes 2} }) = 0.
\]

\ts From the exact sequence (\ref{resolution}) (in the case $\xi\cong\cO$)
we obtain the exact
sequence of cohomology
\[
0 \lra H^0(C^{(g-1)}, T_{C^{(g-1)}} \otimes L) \lra H^1(C, \cO) \otimes
H^0(C^{(g-1)}, L) \lra H^0(C^{(g-1)}, \cI_Z \otimes L^{\otimes 2}) \lra ....
\]
or, by Lemma \ref{Hi},
\begin{equation}
0 \lra H^0(C^{(g-1)}, T_{C^{(g-1)}} \otimes L) \lra H^1(C, \cO)^{\otimes
2} \lra H^0(C^{(g-1)}, \cI_Z \otimes L^{\otimes 2}) \lra ....
\end{equation}
The map $H^1(C, \cO)^{\otimes 2} \lra 
H^0(C^{(g-1)}, \cI_Z \otimes L^{\otimes 2})$
can be described as follows. Choose a basis $\{ D_1, ..., D_g \}$ of $H^1(C,
\cO)$ and think of its elements as translation invariant vector fields on
$Pic^{g-1}C$. Then the isomorphism $H^1(C, \cO) \stackrel{\cong}{\lra}
H^0(C^{(g-1)}, L)$ is defined by $D_i \mapsto \rho^*D_i \theta$ where
$\theta$ is a nonzero element of $H^0(Pic^{g-1}C, \cO_{Pic^{g-1}C}(\T))$ (see
\cite{green84} page 92). Hence the natural map $C^{(g-1)}\ra | L|^*$
is the composition of the rational resolution $\rho$ with the Gauss
map $\T\ra\bP T_0 Pic^{g-1} C\cong\bP^{g-1}$, the map $H^1 (C,\cO )\otimes
L^{-1 }\ra\cI_Z$ is multiplication by global sections of $L$, the map
$H^1(C, \cO)^{\otimes 2} \lra H^0(C^{(g-1)}, \cI_Z \otimes
L^{\otimes 2})$ is defined by $D_i \otimes D_j \mapsto
\rho^*(D_i \theta . D_j \theta)$ and is equal to multiplication
\[
H^0(C^{(g-1)}, L)\otimes H^0(C^{(g-1)}, L)\lra H^0(C^{(g-1)},
L^{\otimes 2})
\]
which factors through $H^0(C^{(g-1)}, \cI_Z \otimes L^{\otimes 2})$.
This multiplication map is induced by the multiplication map
\[
H^0(\bP^{g-1},\cO_{\bP^{ g-1}}(1))\otimes H^0(\bP^{g-1},\cO_{\bP^{
g-1}}(1))\lra H^0(\bP^{g-1},\cO_{\bP^{ g-1}}(2))
\]
via the Gauss map. Since the kernel of this map is $\wedge^2
H^0(\bP^{g-1},\cO_{\bP^{ g-1}}(1))$ and the Gauss map is generically
finite (\cite{sasaki83} page 114), the kernel of $H^0(C^{(g-1)},
L)\otimes H^0(C^{(g-1)}, L)\lra H^0(C^{(g-1)}, L^{\otimes 2})$ is
$\Lambda^2 H^0(C^{(g-1)}, L) =\Lambda^2 H^1 (C ,\cO)$. Hence
$H^0(C^{(g-1)}, T_{C^{(g-1)}} \otimes L)$ is isomorphic to
$\Lambda^2H^1(C, \cO)$.

Now suppose that $\xi^{\otimes 2}\not\cong\cO$. As above, the space
$H^0( C^{( g-1)}, T_{C^{ (g-1) }}\otimes L_{\xi^{\otimes 2}})$ is the
kernel of the map
\[
H^1(C, \cO) \otimes H^0(C^{(g-1)}, L_{\xi^{\otimes 2}}) \lra
H^0(C^{(g-1)}, \cI_Z\otimes L_{\xi }^{\otimes 2})
\]
which sends $D_i\otimes s$ to $\rho^* D_i\theta . s$. Using Lemma
\ref{Hi}, this is immediately seen to be injective.
\qed
\vskip20pt

\subsection{}
To compute the higher cohomology vector spaces of
$T_{ C^{ (g-1)} }\otimes L_{\xi^{\otimes 2} }$, we will first
realize sequence (\ref{resolution}) as a pushforward of an exact
sequence of sheaves on $C^{ (g-1)}\times C$ as follows:
Let $D \subset C^{(g-1)} \times C$ be the universal divisor and let
$p,\;q$ be the projections of $C^{(g-1)}\times C$ on the first and
second factors:
\[
\begin{array}{rccl}
D\hookrightarrow &C^{(g-1)} \times C&\stackrel{q}{\lra}&C\\
&{}^p\downarrow&&\\ &C^{(g-1)}&&
\end{array}
\]
Consider the exact sequence
\begin{equation}
0 \lra p^*L_{\xi^{\otimes 2} } \lra \cO_{C^{(g-1)} \times C}(D)
\otimes p^*L_{\xi^{\otimes 2} } \lra \cO_D(D)\otimes
p^*L_{\xi^{\otimes 2} } \lra 0 .\label{ODL}
\end{equation}
Since $p |_D$ is finite, we have

\subsection{Lemma.}\label{high0} All the higher
direct images $R^i p_* (\cO_D(D) \otimes p^*L_{\xi^{\otimes 2} })$ are
zero for $i>0$.
\vskip20pt

\subsection{}
Also, it is easily seen that 
$
p_*\cO_{C^{(g-1)} \times C}\cong
p_*\cO_{C^{(g-1)} \times C}(D)\cong\cO_{C^{(g-1)}}
$
so that 
\[
p_*
p^*L_{\xi^{\otimes 2} }\cong p_* (\cO_{C^{(g-1)} \times C}(D) \otimes
p^*L_{\xi^{\otimes 2} })\cong L_{\xi^{\otimes 2} }
\]
 by
the projection formula. So we obtain the following exact sequence from
the pushforward by $p$ of sequence (\ref{ODL})
\[
0\lra p_* (\cO_D(D)\otimes p^*L_{\xi^{\otimes 2} } )\lra R^1 p_*
p^*L_{\xi^{\otimes 2} }\lra R^1 p_*
(\cO_{C^{(g-1)} \times C}(D) \otimes p^*L_{\xi^{\otimes 2} } )\lra 0,
\]
or, again by the projection formula,
\[
0\lra p_* \cO_D(D)\otimes L_{\xi^{\otimes 2} } \lra H^1(C,\cO)\otimes
L_{\xi^{\otimes 2} }\lra R^1 p_*\cO_{C^{(g-1)} \times C}(D) \otimes
L_{\xi^{\otimes 2} }\lra 0\; .
\]
By the general theory of Hilbert schemes (\cite{grothendieck60}), we
have
\[
p_* \cO_D(D)\cong T_{C^{(g-1)}},
\]
and the map $T_{C^{(g-1)}}\otimes L_{\xi^{\otimes 2} } \ra
H^1(C,\cO)\otimes L_{\xi^{\otimes 2} }$ obtained by this isomorphism
from the above sequence is equal to the first map in sequence
(\ref{resolution}) because both maps are obtained from the
differential of the morphism $\rho : C^{(g-1)}\ra Pic^{g-1}
C$. Therefore the sequence above is equal to sequence
(\ref{resolution}) and $R^1 p_*\cO_{C^{(g-1)} \times C}(D) =
\cI_Z\otimes L$. So we have proved

\subsection{Lemma.}\label{newres} We have
\[
\begin{array}{ccccccccc}
0 &\lra & p_* (\cO_D(D)\otimes p^*L_{\xi^{\otimes 2} } ) &\lra & R^1
p_* p^*L_{\xi^{\otimes 2} } &\lra & R^1 p_*
(\cO_{C^{(g-1)} \times C}(D) \otimes p^*L_{\xi^{\otimes 2} } ) &\lra & 0 \\
\parallel & &\parallel & &\parallel & &\parallel & &\parallel \\
0 &\lra & T_{C^{(g-1)}}\otimes L_{\xi^{\otimes 2} } &\lra & H^1 (C
,\cO )\otimes L_{\xi^{\otimes 2} } &\lra &\cI_Z\otimes L\otimes
L_{\xi^{\otimes 2} } &\lra & 0
\end{array}
\]
\vskip20pt

From Lemmas \ref{high0}, \ref{newres} and the Leray spectral sequence
we deduce

\subsection{Lemma.}\label{HiTHiD} There are isomorphisms
\[
H^i(C^{(g-1)}, T_{C^{(g-1)}} \otimes L_{\xi^{\otimes 2}
})\stackrel{\cong }{\lra} H^i(D, \cO_D(D)\otimes p^*L_{\xi^{\otimes 2}
})
\]
for all $i\geq 0$.
\vskip20pt

Hence we need to determine the cohomology of $\cO_D (D)\otimes p^*
L_{\xi^{\otimes 2}}$. For this we will consider the push-forward by $q :
C^{(g-1)}\times C\ra C$ of sequence (\ref{ODL}). We first prove

\subsection{Lemma.} For all $i\geq 1$,
\[
R^iq_*(\cO_{C^{(g-1)} \times C}(D) \otimes
p^*L_{\xi^{\otimes 2}}) = 0\; .
\]

\ts It is sufficient to show that, for all $x\in C$ and all $i\geq 1$,
we have $H^i(C^{(g-1)}\times\{ x\}, \cO_{C^{(g-1)}\times\{ x\}}(D)
\otimes p^*L_{\xi^{\otimes 2}}) = 0$. Now $\cO_{C^{(g-1)}\times\{ x\}}(D)
\otimes p^* L_{\xi^{\otimes 2}}\cong \cO_{C^{(g-1)}}(U_x) \otimes
L_{\xi^{\otimes 2}}$ where $U_x\subset C^{(g-1)}$ is the image of
$C^{(g-2)}$ by the morphism $D \mapsto D+x$. We have
$\cO_{C^{(g-1)}}(U_x) \otimes L_{\xi^{\otimes
2}}\cong\cO_{C^{(g-1)}}(U_x) \otimes (L_{\xi^{\otimes 2}}\otimes
L^{-1} )\otimes L$. The invertible sheaf $\cO_{C^{(g-1)}}(U_x) \otimes
(L_{\xi^{\otimes 2}}\otimes L^{-1} )$ is ample because $U_x$ is ample
and $\T_{\xi^{\otimes 2}} -\T$ is algebraically equivalent to
$0$. Since $L$ is the canonical sheaf of $C^{(g-1)}$ (see
\cite{green84} page 88), we have Kodaira vanishing
\begin{equation}\label{HiUx0}
H^i(C^{(g-1)}, \cO_{C^{(g-1)}}(U_x) \otimes L_{\xi^{\otimes 2}}) = 
H^i(C^{(g-1)}, \cO_{C^{(g-1)}}(U_x) \otimes (L_{\xi^{\otimes 2}}
\otimes L^{-1})\otimes L) = 0
\end{equation}
for $i > 0$. \qed
\vskip20pt

Therefore, by the Leray spectral sequence,

\subsection{Lemma.}\label{HCDL} For all $i\geq 0$,
\[
H^i( C^{ (g-1) }\times C, \cO_{C^{(g-1)} \times C}(D) \otimes
p^*L_{\xi^{\otimes 2}})\cong H^i( C, q_*(\cO_{C^{(g-1)}
\times C}(D) \otimes p^*L_{\xi^{\otimes 2}}))\; .
\]

\subsection{}\label{7.16}
Also, pushing (\ref{ODL}) forward by $q$, we obtain the exact sequence
\begin{equation}\label{seqzeta}
0 \lra q_*p^*L_{\xi^{\otimes 2}} \lra q_*(\cO_{C^{(g-1)} \times
C}(D)\otimes p^*L_{\xi^{\otimes 2}})\stackrel{\zeta}{\lra}
q_*(\cO_D(D) \otimes p^*L_{\xi^{\otimes 2}}) \lra
R^1q_*p^*L_{\xi^{\otimes 2}} \lra 0
\end{equation}
and the isomorphisms
\[
R^iq_*(\cO_D(D) \otimes p^*L_{\xi^{\otimes 2}}) \cong
R^{i+1}q_*p^*L_{\xi^{\otimes 2}}
\]
for $i > 0$. We have

\subsection{Lemma.} Suppose that $\xi^{\otimes 2}\cong\cO$. Then the
map $\zeta$ in the exact sequence
\[
0 \lra q_*p^*L \lra q_*(\cO_{C^{(g-1)} \times C}(D)\otimes
p^*L)\stackrel{\zeta}{\lra} q_*(\cO_D(D) \otimes p^*L) \lra
R^1q_*p^*L\lra 0
\]
obtained from sequence (\ref{seqzeta}) is zero.

\ts Let $x$ be an arbitrary element of $C$. Then we have the exact
sequence
\[
0 \lra L \lra \cO_{C^{(g-1)}}(U_x) \otimes L \lra \cO_{U_x}(U_x) \otimes L
\lra 0 \: .
\]
It is sufficient to show that the map on global sections
\[
H^0(C^{(g-1)}, \cO_{C^{(g-1)}}(U_x) \otimes L) \lra H^0(U_x,
\cO_{U_x}(U_x) \otimes L)
\]
obtained from the exact sequence is zero. For this, since
$h^0(C^{(g-1)}, L) = g$ (\ref{Hi}), it is enough to show that
$h^0(C^{(g-1)}, \cO_{C^{(g-1)}}(U_x)
\otimes L) = g$ as well. Since all higher cohomology of
$\cO_{C^{(g-1)}}(U_x) \otimes L$ vanishes (see (\ref{HiUx0})), we need
to show that $\chi(C^{(g-1)}, \cO_{C^{(g-1)}}(U_x) \otimes L) = g$, or
\[
\chi(C^{(g-1)}, L) + \chi(U_x, \cO_{U_x}(U_x) \otimes L) = g \: .
\]
By \cite{green84} page 88, the canonical sheaf of $C^{(g-1)}$ is
$L$. By adjunction, the canonical sheaf of $C^{(g-2)} \cong U_x$ is
$\cO_{U_x}(U_x) \otimes L$. So we need to show
\[
\chi(C^{(g-1)}, \omega_{C^{(g-1)}})+ \chi(C^{(g-2)}, \omega_{C^{(g-2)}})=g
\]
or
\[
h^{g-1, 0}(C^{(g-1)}) - h^{g-1, 1}(C^{(g-1)}) + ... + (-1)^{g-1}h^{g-1,
g-1}(C^{(g-1)}) + \]
\[
+ h^{g-2, 0}(C^{(g-2)}) - h^{g-2, 1}(C^{(g-2)})
 + ... + (-1)^{g-2}h^{g-2, g-2}(C^{(g-2)}) = g \: .
\]
By the weak Lefschetz theorem
\[
h^{g-1,i}(C^{(g-1)}) = h^{g-2,i-1}(C^{(g-2)})
\]
for $i > 1$ because $U_x$ is a smooth and ample divisor in
$C^{(g-1)}$. Therefore we are reduced to showing
\[
h^{g-1, 0}(C^{(g-1)}) - h^{g-1, 1}(C^{(g-1)}) + h^{g-2, 0}(C^{(g-2)}) = g
\: .
\]
We have $H^{g-1, 0}(C^{(g-1)}) \cong H^0(C^{(g-1)},
\omega_{C^{(g-1)}}) \cong H^0(C^{(g-1)}, L) \cong H^1(C, \cO)$ by
Lemma \ref{Hi}. So we are
reduced to showing
\[
h^{g-1, 1}(C^{(g-1)}) = h^{g-2, 0}(C^{(g-2)})
\]
or, by Serre duality,
\[
h^{0,g-2}(C^{(g-1)}) = h^{0, g-2}(C^{(g-2)})
\]
which is true by \cite{macdonald62} (11.1) page 327.
\qed
\vskip20pt

Therefore we have
\subsection{Corollary.} For all $i\geq 0$,
\[
R^iq_*(\cO_D(D) \otimes p^*L)\cong R^{i+1}q_*p^*L
\]
and
\[
q_* p^* L\cong q_*(\cO_{ C^{ (g-1) }\times C} (D)\otimes p^* L)\; .
\]
\vskip20pt

Note that $q_* p^* L\cong H^0 (C^{ (g-1)}, L)\otimes\cO$. Therefore,
combining the above Corollary with Lemmas \ref{Hi} and \ref{HCDL} we obtain

\subsection{Lemma.}\label{cohL} The cohomology groups of the sheaf
$\cO_{C^{(g-1)}\times C}(D)\otimes p^*L$ on $C^{(g-1)}\times C$ are as
follows.
\[
H^i(C^{(g-1)} \times C, \cO_{C^{(g-1)} \times C}(D) \otimes p^*L
)\cong\left\{\begin{array}{lr}
H^1(C,\cO) &\hbox{ if } i= 0 \\
H^1(C, \cO)^{\otimes 2} &\hbox{ if } i= 1 \\
0 &\hbox{ if } i\geq 2
\end{array}\right.
\]

\subsection{}
Now we will compute the cohomology of $\cO_{C^{(g-1)}\times
C}(D)\otimes p^*L_{\xi^{\otimes 2}}$ in the case where $\xi^{\otimes
2}\not\cong\cO$. We have

\subsection{Lemma.}\label{cohLxi} If $\xi^2\not\cong\cO$, the
cohomology groups of the sheaf $\cO_{C^{(g-1)}\times C}(D)\otimes
p^*L_{\xi^{\otimes 2}}$ on $C^{(g-1)}\times C$ are as follows.
\[
H^i(C^{(g-1)} \times C, \cO_{C^{(g-1)} \times C}(D) \otimes p^*
L_{\xi^{\otimes 2}} )\cong\left\{\begin{array}{lr} H^0(\T ,\cO_{\T
}(\T_{\xi^{\otimes 2}})\cong\bC &\hbox{ if } i= 0 \\
\bC^{g^2 - g +1} &\hbox{ if } i= 1 \\
0 &\hbox{ if } i\geq 2
\end{array}\right.
\]

\ts We have $H^0(C^{(g-1)} \times C,\cO_{C^{(g-1)} \times C}(D) \otimes
p^*L_{\xi^{\otimes 2}}) \cong H^0(C^{(g-1)}, p_*(\cO_{C^{(g-1)} \times
C}(D) \otimes p^*L_{\xi^{\otimes 2}}))$ and it can be easily seen that
$L_{\xi^{\otimes 2}} \cong p_*(\cO_{C^{(g-1)} \times C}(D)\otimes
p^*L_{\xi^{\otimes 2}})$. Hence $H^0(C^{(g-1)} \times C,
\cO_{C^{(g-1)} \times C}(D)
\otimes p^*L_{\xi^{\otimes 2}}) \cong H^0(C^{(g-1)}, L_{\xi^{\otimes
2}})$. By Lemma \ref{Hi}, the space $H^0(C^{(g-1)},
L_{\xi^{\otimes 2}})$ is isomorphic to $H^0(\T ,\cO_{\T
}(\T_{\xi^{\otimes 2}})\cong H^0( Pic^{g-1}
C,\T_{\xi^{\otimes 2}} )$ which has dimension $1$.

The fact that $H^i(C^{(g-1)} \times C, \cO_{C^{(g-1)} \times C}
(D)\otimes p^* L_{\xi^{\otimes 2}} ) = 0$ for $i\geq 2$ follows from
Lemma \ref{HCDL}.

To compute the dimension of $H^1(C^{(g-1)} \times C,\cO_{C^{ (g-1)
}\times C}(D) \otimes p^* L_{\xi^{\otimes 2}} )$ we use the following
deformation argument:

Let
\[
\phi : C^{(g-1)}\times C\times Pic^0 C\lra Pic^{g-1} C
\]
be the composition of $\rho\times id\times id$ with the morphism
\[
\T\times C\times Pic^0 C\lra Pic^{g-1} C,\quad (l
,x,\xi)\longmapsto l\otimes\xi\; .
\]
Let $p_{12} : C^{(g-1)}\times C\times Pic^0 C\ra C^{(g-1)}\times C$ be
the projection. Then $p_{12 }^*\cO_{C^{(g-1)}\times
C}(D)\otimes\phi^*\cO_{Pic^{g-1} C} (\T)$ is a family of invertible
sheaves on $C^{(g-1)}\times C$ parametrized by $Pic^0 C$. At $\xi\in
Pic^0 C$, we have
\[
\left( p_{ 12 }^*\cO_{C^{(g-1)}\times C} (D
)\otimes\phi^*\cO_{Pic^{g-1} C} (\T)\right) |_{C^{(g-1)}\times
C\times\{\xi\} }\cong\cO_{C^{(g-1)}\times C}(D)\otimes p^*
L_{\xi^{\otimes 2}}\; .
\]
It follows that, for all $\xi\in Pic^0 C$, we
have $\chi (\cO_{C^{(g-1)}\times C}(D)\otimes p^* L_{\xi^{\otimes 2}})
=\chi (\cO_{C^{(g-1)}\times C}(D)\otimes p^* L )$.  Hence $\chi
(\cO_{C^{(g-1)}\times C}(D)\otimes p^* L_{\xi^{\otimes 2}}) = g - g^2$
by Lemma \ref{cohL}. Therefore $h^0 ( C^{ (g-1) }\times C,
\cO_{C^{(g-1)}\times C}(D)\otimes p^*L_{\xi^{\otimes 2}}) - h^1 (C^{
(g-1) }\times C, \cO_{C^{(g-1)}\times C}(D)\otimes p^*L_{\xi^{\otimes
2}}) = g - g^2$ and $h^1 ( C^{ (g-1) }\times C,\cO_{C^{(g-1)}\times
C}(D)\otimes p^*L_{\xi^{\otimes 2}}) = g^2 - g + 1$.
\qed
\vskip20pt

We are now ready to compute the cohomology of the sheaf
$\cO_D (D)\otimes p^* L_{\xi^{\otimes 2}}$. We first consider the case
where $\xi^{\otimes 2}\not\cong\cO$. From Lemma \ref{Hi} it follows

\subsection{Lemma.} If $\xi^{\otimes 2}\not\cong\cO$, then
\[
R^iq_*p^*L_{\xi^{\otimes 2}} = 0
\]
for $i\geq 1$.

\subsection{}
So, if $\xi^{\otimes 2}\not\cong\cO$, from \ref{7.16} we obtain the exact sequence
\begin{equation}\label{exseqq}
0 \lra q_*p^*L_{\xi^{\otimes 2}} \lra q_*(\cO_{C^{(g-1)} \times C}(D)
\otimes p^*L_{\xi^{\otimes 2}})
\stackrel{\zeta}{\lra}
q_*(\cO_D(D) \otimes p^*L_{\xi^{\otimes 2}}) \lra 0
\end{equation}
and
\[
R^iq_*(\cO_D(D) \otimes p^*L_{\xi^{\otimes 2}}) = 0
\]
for $i > 0$.

Therefore, again by the Leray spectral sequence,

\subsection{Lemma.}\label{cohLDq} If $\xi^{\otimes 2}\not\cong\cO$,
then, for $i\geq 0$,
\[
H^i( C^{ (g-1) }\times C, p^*L_{\xi^{\otimes 2}})\cong H^i( C, q_*
p^*L_{\xi^{\otimes 2}})
\]
and
\[
H^i( D, \cO_D(D) \otimes p^*L_{\xi^{\otimes 2}})\cong H^i (C,
q_*(\cO_D(D) \otimes p^*L_{\xi^{\otimes 2}}))\; .
\]
In particular, $H^i( D, \cO_D(D) \otimes
p^* L_{\xi^{\otimes 2}}) = H^i(C^{ (g-1) }\times C, p^*L_{\xi^{\otimes 2}})=
0$ for $i\geq 2$.

\subsection{}
Hence the exact sequence (\ref{exseqq}) gives the exact sequence of cohomology
\[
0\lra H^0 (C, q_* p^* L_{\xi^{\otimes 2}})\lra H^0 (C, q_* (\cO_{C^{
(g-1) }\times C} (D)\otimes p^* L_{\xi^{\otimes 2}})) \lra H^0 (C, q_*
(\cO_D (D)\otimes p^* L_{\xi^{\otimes 2}})) \lra\phantom{0\; .}
\]
\[
\phantom{0}\lra H^1 (C, q_* p^*
L_{\xi^{\otimes 2}})\lra H^1 (C, q_* (\cO_{C^{ (g-1) }\times C}
(D)\otimes p^* L_{\xi^{\otimes 2}})) \lra H^1 (C, q_* (\cO_D
(D)\otimes p^* L_{\xi^{\otimes 2}})) \lra 0\; .
\]
By Lemmas \ref{H0TL}, \ref{HiTHiD} and \ref{cohLDq}, we have $H^0 (C, q_*
(\cO_D (D)\otimes p^* L_{\xi^{\otimes 2}})) = H^0 (D,\cO_D (D)\otimes
p^* L_{\xi^{\otimes 2}}) = 0$. Therefore, by Lemmas \ref{HCDL},
and \ref{cohLDq},
\[
0\lra H^1 (C^{ (g-1) }\times C, p^*
L_{\xi^{\otimes 2}})\lra H^1 (C^{ (g-1) }\times C, \cO_{C^{ (g-1) }\times C}
(D)\otimes p^* L_{\xi^{\otimes 2}}) \lra 
\]
\[\lra H^1 (C^{ (g-1) }\times C, \cO_D
(D)\otimes p^* L_{\xi^{\otimes 2}}) \lra 0\; .
\]
Therefore, from the above Lemma, we can deduce

\subsection{Lemma.}\label{h1D} We have
\[
h^1 (D,\cO_D (D)\otimes p^* L_{\xi^{\otimes 2}}) = (g-1)^2\; .
\]

\subsection{Proof of part 1 of Theorem \ref{thmsing}.} The exact
sequence in part 1 of Theorem \ref{thmsing} is obtained from the long
exact sequence of cohomology of sequence (\ref{resolution}) by using
Lemmas \ref{Hi} and \ref{H0TL}. Then it follows from Lemmas
\ref{HiTHiD} and \ref{h1D} that $h^1 ( C^{ (g-1) }, T_{C^{ (g-1)
}}\otimes L_{\xi^{\otimes 2}}) = (g-1)^2$.
\qed

\subsection{}
From now on we will suppose that $\xi^{\otimes 2}\cong\cO$. We will
compute the higher cohomology of $\cO_D (D)\otimes L$:

\subsection{Lemma.}\label{lemT}
There are isomorphisms:
\[
H^i(D, \cO_D(D)\otimes p^*L)\stackrel{\cong}{\lra} 
H^{i+1}(C^{(g-1)} \times C, p^*L)
\]
for all $i\geq 1$.

\ts Consider the long exact sequence of cohomology of (\ref{ODL}) and
use Lemma \ref{cohL} to obtain the exact sequence
\begin{eqnarray}
0 \lra H^0(C^{(g-1)} \times C, p^* L) \lra
H^1 (C,\cO)\phantom{^{\otimes 2}} \lra H^0 ( D,\cO_D (D)\otimes
p^*L)\lra\phantom{0}\nonumber \\
\label{exseq5}\phantom{0}\lra H^1 ( C^{(
g-1)}\times C, p^*L) \lra H^1 (C ,\cO )^{\otimes 2} \lra H^1 ( D,\cO_D
(D)\otimes p^*L)\lra \phantom{0}\\
\phantom{0} \lra H^2 ( C^{(
g-1)}\times C, p^*L) \lra 0\phantom{
H^1 (C,\cO)\phantom{^{\otimes 2}} \lra H^0 ( D,\cO_D (D)\otimes
p^*L)\lra 0}\nonumber
\end{eqnarray}
and the isomorphisms
\[
H^i(D, \cO_D(D) \otimes p^*L) \stackrel{\cong}{\lra} 
H^{i+1}(C^{(g-1)} \times C, p^*L)
\]
for all $i\geq 2$.

It remains to do the case $i=1$. By Lemmas \ref{H0TL} and \ref{HiTHiD}
we have:
\[
H^0(D, \cO_{D}(D) \otimes p^*L) \cong H^0(C^{(g-1)}, T_{C^{(g-1)}}
\otimes L)\cong \Lambda^2H^1(C, \cO)\; .
\]
By the K\"unneth isomorphism and Lemma \ref{Hi},
\[
H^0( C^{ (g-1) }\times C, p^* L)\cong H^0 (C^{ (g-1) }, L)\cong
H^1(C,\cO)
\]
and
\[ \begin{array}{ccl}
H^1(C^{(g-1)} \times C, p^*L) & \cong & \left( H^0(C^{(g-1)}, L) \otimes
H^1(C, \cO) \right) \oplus H^1(C^{(g-1)}, L) \\
 & \cong & H^1(C, \cO)^{\otimes 2}\oplus \Lambda^2H^1(C, \cO).
\end{array}
\]
Putting these results in sequence (\ref{exseq5}) we obtain
\begin{eqnarray}
0 \lra \Lambda^{2}H^1(C, \cO) \lra H^1(C, \cO)^{\otimes
2}\oplus\Lambda^{ 2} H^1( C, \cO) \lra H^1(C, \cO)^{\otimes
2}\lra\phantom{ 0\: .} \\ \lra H^1 (D,\cO_D(D) \otimes p^*L)
\stackrel{\psi}{\lra} H^2(C^{(g-1)} \times C, p^*L)\lra 0 \:
. \nonumber
\end{eqnarray}
and it follows that $\psi$ is an isomorphism. This concludes the proof of the 
Lemma. \qed

\subsection{Proof of part 2 of Theorem \ref{thmsing}.} We take the
cohomology of the exact sequence (\ref{resolution}):
\begin{eqnarray}
0\longrightarrow H^0(T_{C^{(g-1)}}\otimes L)\longrightarrow
H^1(C,\cO)\otimes H^0(C^{(g-1)}, L)\lra
H^0(C^{(g-1)},\cI_Z\otimes L^{\otimes 2}) \lra
\label{eq45}\\
\phantom{0}\lra H^1(T_{C^{(g-1)}}\otimes L)\lra
H^1(C,\cO)\otimes H^1 (C^{ (g-1)},L )\lra\ldots
\phantom{H^0(C^{(g-1)},\cI_Z\otimes L^{\otimes 2})\;\;} .\nonumber
\end{eqnarray}

\subsection{Lemma.}\label{surj}
The map:
\[
H^1(T_{C^{g-1}}\otimes L)\longrightarrow
H^1(C,\cO)\otimes H^1(C^{g-1},L)
\]
in sequence (\ref{eq45}) is surjective. 

\ts Since the fibers of $p$ are one
dimensional, $R^2p_*=0$. The Leray spectral sequence then shows that we
have a surjective map:
\[
H^2(C^{(g-1)} \times C, p^*L) \surj
H^1(C^{(g-1)}, R^1p_*(p^*L))\cong 
H^1(C^{(g-1)},(R^1p_*\cO)\otimes L)
\cong H^1(C,\cO)\otimes H^1(C^{(g-1)},L),
\]
where the isomorphisms are given by the projection formula and the
fact that $R^1p_*\cO\cong \cO\otimes H^1(C,\cO)$. By Lemmas
\ref{HiTHiD} and \ref{lemT} we have an isomorphism
$H^1(T_{C^{(g-1)}}\otimes L)\cong H^2(C^{(g-1)}
\times C, p^*L)$ and it is easily seen, using Lemma \ref{newres}, that
the composition 
$H^1(T_{C^{(g-1)}}\otimes L)\cong H^2(C^{(g-1)} \times C, p^*L)
\surj H^1(C,\cO)\otimes
H^1(C^{(g-1)},L)$ 
coincides with the map in sequence (\ref{eq45}). This
concludes the proof.
\qed

\subsection{}
Using Lemmas \ref{Hi}, \ref{H0TL}, \ref{HiTHiD} and Lemma \ref{lemT}
together with the K\"unneth isomorphism, we obtain the following exact
sequence from sequence (\ref{eq45})
\[ \begin{array}{l}
0 \lra \Lambda^2H^1(C, \cO) \lra H^1(C, \cO)^{\otimes 2} \lra
H^0(C^{(g-1)}, \cI_Z \otimes L^{\otimes 2}) \lra \\
\phantom{0}\lra 
\left( \Lambda^2H^1(C, \cO) \otimes H^1(C, \cO) \right) \oplus
\Lambda^{3}H^1(C, \cO)
\lra \Lambda^2H^1(C, \cO) \otimes H^1(C, \cO) \lra 0
\end{array} 
\]
and Theorem \ref{thmsing} follows.
\qed


\providecommand{\bysame}{\leavevmode\hbox to3em{\hrulefill}\thinspace}

\end{document}